\documentclass[11pt]{amsart}
\usepackage{amsmath,amssymb,mathrsfs, xcolor, geometry}
\usepackage{hyperref}
\hypersetup{linkbordercolor = {white}, citebordercolor = {white},}
\usepackage[alphabetic]{amsrefs}

\usepackage[alphabetic]{amsrefs}

\numberwithin{equation}{section}

\theoremstyle{definition}
\newtheorem{thm}[equation]{Theorem} 
\theoremstyle{definition}
\newtheorem{definition}[equation]{Definition}
\newtheorem{prop}[equation]{Proposition} 
\newtheorem{cor}[equation]{Corollary}

\newtheorem{lemma}[equation]{Lemma}

\newcommand{\ovl}[1]{\overline{#1}}

\newcommand{\pair}[1]{\langle #1\rangle}
\newcommand{\set}[1]{\left\{#1\right\}}
\newcommand{\ppair}[1]{\left\langle #1\right\rangle}
\newcommand{\pr}[1]{\left(#1\right)}
\newcommand{\abs}[1]{\left|#1\right|}
\newcommand{\bb}[1]{\mathbb{#1}} \newcommand{\td}[1]{\widetilde{#1}}

\newcommand{\tr}{\text{tr}}
\newcommand{\D}{\Delta}

\newcommand{\n}{\nabla}

\newcommand{\sbst}{\subseteq}
\newcommand{\h}{\textbf H}

\newcommand{\bd}{\partial}

\newcommand{\N}{\mathbf{n}}
\newcommand{\cH}{\mathcal{H}}

\newcommand{\cC}{\mathcal C}

\newcommand{\cS}{\mathcal S}

\usepackage{enumitem}
\setlist[enumerate]{itemsep=3pt,topsep=5pt}
\allowdisplaybreaks[2]
\newcommand{\R}{\mathbb R}
\newcommand{\C}{\mathbb C}

\newcommand{\op}{\mathrm{op}}
\newcommand{\eps}{\varepsilon}
\newcommand{\norm}[1]{\left\lVert#1\right\rVert}
\newcommand{\ip}[2]{\left\langle#1,#2\right\rangle}

\DeclareMathOperator{\diag}{diag}
\DeclareMathOperator{\RePart}{Re}
\DeclareMathOperator{\Pf}{Pf}

\title{Almost isoclinic Lagrangian mean curvature flow 
}

\author{Tang-Kai Lee}
\author{Mu-Tao Wang}

\address{Department of Mathematics, Columbia University, New York, NY 10027, USA}
\email{leetk@math.columbia.edu and mtwang@math.columbia.edu}

\date{\today}

\begin{document}
\begin{abstract}
We study the mean curvature flow of almost isoclinic Lagrangian submanifolds, a class introduced in \cite{LW26}. 
The almost isoclinic condition is formulated in terms of the Gauss map into the Lagrangian Grassmannian. This class includes spiral-shaped submanifolds that are not graphical over any fixed Lagrangian \(n\)-plane in \(\mathbb{C}^n\). 
Consequently, their mean curvature flow cannot be described globally by the parabolic special Lagrangian equation studied in previous work. 
We prove that, under suitable asymptotic conditions, these submanifolds unwind and flatten along the mean curvature flow. 
This gives a higher-dimensional analogue of the unwinding behavior of planar curves under classical curve shortening flow and leads to topological consequences.
\end{abstract}
\maketitle

\section{\bf Introduction}
We study the mean curvature flow of almost isoclinic Lagrangian submanifolds in $\mathbb C^n$.
The almost isoclinic condition for Lagrangian submanifolds and the graphical condition for hypersurfaces both constrain the Gauss map to a maximal convex region of the corresponding Grassmannian. 
For hypersurfaces, this region is an open hemisphere of the sphere; for Lagrangian submanifolds, it is the almost isoclinic region of the Lagrangian Grassmannian. 
Convexity assumptions for Lagrangian graphs have traditionally been formulated in terms of their potentials. 
The almost isoclinic condition is formulated directly through the Gauss map and allows submanifolds that are not graphical over any fixed Lagrangian plane; see, for example, \cite{LW26}*{Section~4.2} and \cite{LW26}*{Figure~3}.
This formulation allows us to apply the theorem in \cite{W03} that the Gauss map evolves by harmonic map heat flow along mean curvature flow to show that the almost isoclinic condition is preserved. 
In this article, we show that this condition also yields the estimates for the second fundamental form needed to establish long-time existence and convergence of the flow.

Existing long-time existence and convergence results for Lagrangian mean curvature flow generally fall into three categories: (1) graphical settings \cite{SW02, CCH12, CCY13, TTW24, TTW25}, (2) symmetry assumptions that reduce to a curve flow \cite{S20, LO20}, and (3) results under suitable smallness assumptions \cite{L12, LT24, CJL}. Here we leave aside curve shortening flow on surfaces, which is a special case of Lagrangian mean curvature flow but is usually studied separately \cite{GH86, G89, G90, CZ98}. 
Our results require only the general open condition given by a lower bound on the isoclinicity function $\Lambda$, without any graphicality assumption.

In this paper, we use the following conventions unless otherwise stated.
All submanifolds are assumed to be smoothly immersed, complete, and connected in $\bb R^{2n}$.
In particular, we do not a priori assume embeddedness.
A special Lagrangian cone means a finite union of minimal Lagrangian cones with the same Lagrangian angles.
A Lagrangian submanifold $M$ is called an asymptotically conical Lagrangian if it is asymptotic to a special Lagrangian regular cone.

The first theorem concerns the case where the ambient space is $\mathbb{C}^n$ with $n$ even.
In this case, it suffices to assume a positive lower bound on the isoclinicity function $\Lambda.$

\begin{thm}
\label{thm:even-flow}
    Let $n\ge 2$ be even and $M$ be an $n$-dimensional asymptotically conical almost isoclinic Lagrangian in $\bb C^n.$
    If there exists $c>0$ such that $\Lambda\ge c>0,$ then there exists a unique Lagrangian mean curvature flow $M_t$ starting from $M$ that exists for all $t\ge 0.$
    Moreover, a rescaling of the flow converges smoothly and subsequentially to a Lagrangian $n$-plane.
\end{thm}

The second result holds for almost isoclinic Lagrangians with $\Lambda\ge 1/2$ in any dimension $n$.
We state it only for odd $n$, since the even case is already covered by the previous theorem.

\begin{thm}
\label{thm:odd-flow}
    Let $n\ge 3$ be odd and $M$ be an $n$-dimensional asymptotically conical, zero-Maslov, almost isoclinic Lagrangian in $\bb C^n.$
    If $\Lambda\ge 1/2,$ then there exists a unique Lagrangian mean curvature flow $M_t$ starting from $M$ that exists for all $t\ge 0.$
    Moreover, a rescaling of the flow converges smoothly and subsequentially to a Lagrangian $n$-plane.
\end{thm}

Finally, assuming an $SO(n)$-symmetry, we analyze the underlying modified curve shortening flow and prove long-time and convergence for the corresponding Lagrangian mean curvature flows.

\begin{thm}
\label{thm:symmetric-flow}
    Let $n\ge 2$ and $M$ be an $n$-dimensional asymptotically conical, $SO(n)$-symmetric, almost isoclinic Lagrangian in $\bb C^n.$
    If there exists $c>0$ such that $\Lambda\ge c>0,$ then there exists a unique Lagrangian mean curvature flow $M_t$ starting from $M$ that exists for all $t\ge 0.$
    Moreover, a rescaling of the flow converges smoothly and subsequentially to a Lagrangian $n$-plane.
\end{thm}

We remark that neither embeddedness nor a specific topological type is assumed in any of the theorems above.
In fact, the proofs of Theorems~\ref{thm:even-flow} and~\ref{thm:symmetric-flow} show that, under the almost isoclinicity assumption, the initial Lagrangian is already diffeomorphic to $\bb R^n$. 
By contrast, for odd-dimensional almost isoclinic Lagrangians, this topological conclusion can be established only after long-time existence of the flow has been proved.
For this reason, we record the following topological consequence.

\begin{cor}
    Suppose a submanifold $M$ satisfies the assumptions in one of Theorems~\ref{thm:even-flow}, \ref{thm:odd-flow}, and~\ref{thm:symmetric-flow}.
    Then $M$ is diffeomorphic to $\bb R^n.$
\end{cor}

We briefly describe the proof ideas and the structure of the paper.
Each of the main theorems follows from three steps.
We will first prove that the mean curvature flow starting from a given almost isoclinic Lagrangian in the theorems exists for all positive  time.
Next, we prove that the long-time existence of the flow forces the Lagrangian becomes graphical for large enough time, which is done in Section~\ref{sec:G-criterion}.
This then reduces the task to the case of a graphical Lagrangian mean curvature flow, which is dealt with in Section~\ref{sec:G-MCF}.

For each of the three main theorems, the long-time existence is established by different approaches.
In Section~\ref{sec:even-G}, we show that for a Lagrangian of even dimension $n=2m$ in $\mathbb{C}^n$, almost isoclinicity implies graphicality over a complex subspace of real dimension $n$ (hence of complex dimension $m$). Long-time existence then follows, in Corollary~\ref{cor:tangent-flow-structure}, from a blow-up argument combined with a structure theorem for tangent flows. For an odd-dimensional Lagrangian, the long-time existence is obtained by a pointwise curvature estimate in Section~\ref{sec:odd}.
Finally, the $SO(n)$-symmetric case is based on a non-collapsing estimate for the underlying modified curve shortening flow, which is in Section~\ref{sec:sym}.

\subsection*{Acknowledgment}
T.-K. Lee acknowledges the support from the AMS-Simons Travel Grant.
M.-T.~Wang is supported in part by the National Science Foundation under Grant DMS-2404945. This work was supported by grants from the Simons Foundation [MPS-TSM-00007411, M.-T.~W.].

\subsection*{Disclosure of AI use} M.-T.~Wang used ChatGPT 6.1 Sol to assist with literature searches, explore ideas, and perform calculations during the preparation of this work. For Sections 3.2 and 4.1, M.-T.~Wang developed the ideas, verified the theorem statements, and prepared drafts, which were then typeset with AI assistance. 
After the authors had obtained the results, ChatGPT was used to check the validity of some arguments.
All proofs and calculations were independently verified by the authors.

\section{\bf Preliminary}

\subsection{Almost isoclinic condition}

Lee--Wang \cite{LW26} introduced the isoclinicity function $\Lambda$ and the almost isoclinic condition for an oriented Lagrangian plane $L$ in $\mathbb{C}^n$,
with respect to a fixed oriented base plane $\mathbb R^n\subset\mathbb C^n$ for $n\geq 2$. 
The almost isoclinic condition uses a phase rotation $e^{-i\varphi}L$ that is graphical over $\mathbb R^n$ and imposes inequalities on the eigenvalues of the associated symmetric graph matrix, together with an orientation condition on the projection.  
In this section, we construct an angle basis for each $L$ and express the isoclinicity function $\Lambda$ and the almost isoclinic condition in terms of the corresponding angles of $L$. 
We describe the ambiguity in these angles and verify that the resulting conditions and quantities are well-defined.

\subsubsection{Definitions using the graph eigenvalues}

As in \cite{LW26}, we consider an oriented Lagrangian plane $L$ in $\mathbb{C}^n$ as an element of the Lagrangian Grassmannian $\operatorname{Lag}^+(n) \simeq {U}(n)/{SO}(n)$ and represent $L$ by $X+iY\in U(n)$, whose columns form an oriented orthonormal basis for $L$. 
Here $X$ and $Y$ are real matrices.
Changing this basis replaces $X+iY$ by $(X+iY)Q$ for some $Q\in SO(n)$. 

\begin{definition}[See Section~2.3 of \cite{LW26}]
\label{def:AI}
Let $L\in {\rm Lag}^+(n)$ and represent $L$ by $X+iY$ as above.
\begin{enumerate}
\item The \emph{isoclinicity function} $\Lambda$ on $\operatorname{Lag}^+(n)$ is defined by
\begin{equation*}
\Lambda(L)=\det\bigl(\wedge^2X+\wedge^2Y\bigr).
\end{equation*} 
This is independent of the chosen oriented orthonormal basis. 
Moreover, writing 
$X_\varphi+iY_\varphi=e^{-i\varphi}(X+iY)$, we have 
\[
\wedge^2X_\varphi+\wedge^2Y_\varphi
=\wedge^2X+\wedge^2Y.
\] 
and thus $\Lambda(e^{-i\varphi}L)=\Lambda(L)$ for every $\varphi\in \mathbb R$.

\item The plane $L$ is \emph{almost isoclinic}, written $L\in\mathcal{AI}(n)$, if there exists $\varphi\in \mathbb R$ such that 
\[
\det X_\varphi>0
\quad\text{and}\quad
1+\lambda_j\lambda_k>0\qquad\text{for every }j<k,
\]
where $\lambda_1,\cdots, \lambda_n$ are the eigenvalues of $Y_\varphi X_\varphi^{-1}$. 
In this case, 
\begin{equation*}
\Lambda(L)=
\prod_{i<j}
\frac{1+\lambda_i\lambda_j}
{\sqrt{1+\lambda_i^2}\sqrt{1+\lambda_j^2}}>0.
\end{equation*}
\end{enumerate}
\end{definition}

\subsubsection{An angle basis for an arbitrary oriented Lagrangian plane}

The following lemma amounts to writing a representative of $L$ in the form
\[
X+iY
=\pr{u_1~\cdots~u_n}\operatorname{diag}\bigl(e^{i\theta_1},\ldots,e^{i\theta_n}\bigr).
\]

\begin{lemma}\label{lem:angle-basis}
For any oriented Lagrangian plane \(L\subset\mathbb C^n\), there exist a positively oriented orthonormal basis \(e_1,\ldots,e_n\) of \(L\), a positively oriented orthonormal basis \(u_1,\ldots,u_n\) of the standard \(\mathbb R^n\), and real numbers \(\theta_1,\ldots,\theta_n\) such that for any $j=1,\ldots,n,$
\begin{equation}\label{eq:adapted-frame}
e_j=e^{i\theta_j}u_j.
\end{equation}
\end{lemma}

\begin{proof}
Choose a positively oriented orthonormal basis of \(L\), with column matrix \(X+iY\). 
Since \(L\) is Lagrangian, this matrix is unitary, so \(\det(X+iY)\ne0\). 
Thus, the real polynomial \(\det(X+tY)\) is not identically zero. 
Choose \(t\in\mathbb R\) for which this determinant is nonzero, and choose \(\phi\in(-\pi/2,\pi/2)\) with \(\tan\phi=t\).
The real projection of \(e^{-i\phi}L\), whose matrix is \(\cos\phi\,X+\sin\phi\,Y\), is then invertible.

Consequently, \(e^{-i\phi}L\) is the graph
\(\{x+iSx:x\in\mathbb R^n\}\) of a real linear map \(S\).
The Lagrangian condition gives
\[
0=\omega(x+iSx,y+iSy)
=\langle x,Sy\rangle-\langle Sx,y\rangle,
\]
so \(S\) is symmetric.
Choose a positively oriented orthonormal
eigenbasis \(u_1,\ldots,u_n\), with \(Su_j=\lambda_j u_j\), and set \(\alpha_j=\arctan\lambda_j\). 
Then
\[
f_j=\frac{u_j+iSu_j}{\sqrt{1+\lambda_j^2}}
=e^{i\alpha_j}u_j
\]
form an orthonormal basis of \(e^{-i\phi}L\).
If this basis has the opposite orientation from the one transported from \(L\), replace \(f_1\) by \(-f_1\) and \(\alpha_1\) by \(\alpha_1+\pi\), keeping all \(u_j\) fixed. 
Finally, set \(e_j=e^{i\phi}f_j\) and \(\theta_j=\phi+\alpha_j\).
Both bases now have the required orientations.
\end{proof}

\subsubsection{ Isoclinicity function \texorpdfstring{$\Lambda$}{Lambda} and almost isoclinicity in terms of the angles}

\begin{prop}
\label{prop:angle-characterization}
For every choice of angles and bases in Lemma~\ref{lem:angle-basis},
\begin{equation}\label{eq:lambda-angles}
\Lambda(L)=\prod_{j<k}\cos(\theta_j-\theta_k).
\end{equation}
Moreover, the plane $L$ is almost isoclinic if and only if there is a choice of the angles and bases in Lemma~\ref{lem:angle-basis} such that for every $j$ and $k,$
\begin{equation}\label{eq:positive-cosines}
\cos(\theta_j-\theta_k)>0.
\end{equation}

\end{prop}

\begin{proof}
To prove \eqref{eq:lambda-angles} for an arbitrary choice of the oriented angle bases, set $Q=\pr{u_1~\cdots~u_n}\in SO(n)$.
Then $X=Q\,\operatorname{diag}(\cos\theta_j)$ and $Y=Q\,\operatorname{diag}(\sin\theta_j).$
Consequently, relative to the standard basis $e_j\wedge e_k$,
$j<k$, of $\wedge^2\mathbb R^n$,
\[
\wedge^2X+\wedge^2Y
=(\wedge^2Q)\,
\operatorname{diag}_{j<k}
\bigl(\cos\theta_j\cos\theta_k+\sin\theta_j\sin\theta_k\bigr).
\]
Since $\det(\wedge^2Q)=1$, taking determinants gives \eqref{eq:lambda-angles}.
This also proves that the product is independent of all the choices.

Now suppose $L$ is almost isoclinic. 
Apply the construction in the proof of Lemma~\ref{lem:angle-basis} with
$U=X_\varphi+iY_\varphi,$ where
\[
X_\varphi^{-1}Y_\varphi v_j=\lambda_jv_j
\quad\text{ and }\quad
\det(v_1,\dots,v_n)=1.
\]
The eigenvalues here agree with those of $Y_\varphi X_\varphi^{-1}$.
The construction gives
\[
u_j=\sqrt{1+\lambda_j^2}\,X_\varphi v_j
\quad\text{ and }\quad
\theta_j=\varphi+\arctan\lambda_j.
\]
The real vectors $u_j$ are orthonormal, and now their orientation is already positive since
\[
\det\pr{u_1~\cdots~u_n}
=\det X_\varphi\,\prod_{j=1}^n\sqrt{1+\lambda_j^2}>0.
\]
 Both bases have the required orientations, and
\[
\cos(\theta_j-\theta_k)
=\frac{1+\lambda_j\lambda_k}
{\sqrt{(1+\lambda_j^2)(1+\lambda_k^2)}}>0.
\]

Conversely, suppose the oriented bases in Lemma~\ref{lem:angle-basis} have been chosen so that \eqref{eq:positive-cosines} holds.
Set
\[
Q=\pr{u_1~\cdots~u_n}\in SO(n),
\quad
\varphi=\theta_1,
\quad\text{and}\quad
\alpha_j=\theta_j-\varphi.
\]
Then $\cos\alpha_j>0$ for every $j$. 
Using $e^{i\theta_j}u_j$ as the columns of the representative of $L$ gives $X_\varphi=Q\,\operatorname{diag}(\cos\alpha_j)$ and $Y_\varphi=Q\,\operatorname{diag}(\sin\alpha_j).$
It follows that
\[
\det X_\varphi=\prod_j\cos\alpha_j>0,
\quad\text{ and }\quad
Y_\varphi X_\varphi^{-1}
=Q\,\operatorname{diag}(\tan\alpha_j)\,Q^T.
\]
Its eigenvalues $\lambda_j=\tan\alpha_j$ satisfy
\[
1+\lambda_j\lambda_k
=\frac{\cos(\theta_j-\theta_k)}
{\cos\alpha_j\cos\alpha_k}>0.
\]
Thus, $L$ is almost isoclinic.
\end{proof}

We prove the following lemma that shows the angles can be chosen in an interval of length less than \texorpdfstring{$\pi/2$}{pi/2} for an almost isoclinic $L$. 

\begin{lemma}\label{lem:short-interval}
Let $L$ be an almost isoclinic oriented Lagrangian plane.
There are angles and oriented bases as in Lemma~\ref{lem:angle-basis} such that
\begin{equation}\label{eq:angle-spread}
\max_j\theta_j-\min_j\theta_j<\frac{\pi}{2}.
\end{equation}

\end{lemma}

\begin{proof}
By Proposition~\ref{prop:angle-characterization}, choose angles satisfying $\cos(\theta_j-\theta_k)>0$ for every $j$ and $k.$
In particular, $\cos(\theta_j-\theta_1)>0$.
For each $j$, choose an integer $k_j$ such that
\[
\widetilde\theta_j=\theta_j+2\pi k_j
\quad\text{satisfies}\quad
\theta_1-\frac{\pi}{2}
<\widetilde\theta_j
<\theta_1+\frac{\pi}{2}.
\]
We may take $k_1=0$. 
Such a choice exists because $\cos(\theta_j-\theta_1)>0$. 
Moreover,
\[
e^{i\widetilde\theta_j}u_j=e^{i\theta_j}u_j,
\quad\text{ and }\quad
\cos(\widetilde\theta_i-\widetilde\theta_j)
=\cos(\theta_i-\theta_j)>0.
\]

Set
\[
a=\min_j\widetilde\theta_j
\quad\text{ and }\quad
b=\max_j\widetilde\theta_j.
\]
The preceding bounds imply $0\leq b-a<\pi$.
Since $a$ and $b$ are two of the adjusted angles, pairwise positivity also gives $\cos(b-a)>0$.
On the interval $[0,\pi)$, cosine is positive exactly on $[0,\pi/2)$. 
Therefore, $b-a<\pi/2.$
Relabeling the adjusted angles as $\theta_j$ proves \eqref{eq:angle-spread}.
\end{proof}

\subsection{Asymptotically conical submanifolds}
\label{sec:AC}

In this section, we explain the precise definition of asymptotic conicality that is used in this paper.
We will explain that asymptotically conical submanifolds defined in this way are exactly the ones that satisfy the definition used in Su's short-time existence result \cite{S20}.
Recall that all submanifolds are assumed to be smooth, complete, and connected unless otherwise stated.
We will deal with general regular cones in this section and restrict to the (special) Lagrangian case later.

\begin{definition}
	\label{def:AC}
	Let $\cC$ be a regular cone in $\bb R^N,$ in the sense that $\cC\cap \bb S^{N-1}$ is a smooth compact submanifold of $\bb S^{N-1}.$
	We say a submanifold $F\colon M\to\bb R^N$ is asymptotically conical with rate $\beta<2$ and cone $\cC$ if there exist a compact set $K\sbst M,$ a positive number $R,$ and a smooth section $u\in C^\infty\pr{\cC\setminus \ovl B_R; N\cC}$ such that $F\pr{M\setminus K}
	= \set{x+u(x): x\in\cC\setminus\ovl B_R}$
	and for any $k\ge 0,$
	\begin{align}\label{eq:AC-rate}
	|\n^k_{\cC}u|=O\pr{|x|^{\beta-1-k}}
	\end{align}
	as $|x|\to\infty.$
\end{definition}

We remark that such an asymptotically conical submanifold is always a proper immersion.

\begin{lemma}
	A submanifold $M$ is asymptotically conical with rate $\beta<2$ and cone $\cC$ in the sense of Definition~\ref{def:AC} if and only if it is asymptotically conical with rate $\beta<2$ and cone $\cC$ in the sense of \cite{S20}*{Definition~2.2}.
\end{lemma}

\begin{proof}
	Let $\Sigma:=\cC\cap \bb S^{N-1}$ be the link of $\cC.$
	
	Suppose $M$ is asymptotically conical with rate $\beta<2$ and cone $\cC$ in the sense of Definition~\ref{def:AC}.
	Then fixing a diffeomorphism $\varphi\colon (R,\infty)\times\Sigma\to M\setminus K$ by the map $u,$ i.e., defining $\varphi(x):=x+u(x),$ shows that \eqref{eq:AC-rate} is equivalent to \cite{S20}*{Equation (1)}, which says
	\begin{align}\label{eq:Su-AC}
	|\n^k(F\circ\varphi - \iota_{\cC})|
	= O\pr{r^{\beta-1-k}}
	\end{align}
	as $r\to\infty,$ where $\iota_\cC\colon (0,\infty)\times\Sigma\to\bb R^N$ parametrizes $\cC$ in the standard way.
	
	Suppose $M$ is asymptotically conical with rate $\beta<2$ and cone $\cC$ in the sense of \cite{S20}*{Definition~2.2}.
	This means that there exist $R>0,$ a compact set $K\sbst M,$ and a diffeomorphism $\varphi\colon (R,\infty)\times\Sigma\to M\setminus K$ such that \eqref{eq:Su-AC} holds.
	Consider the tubular neighborhood map $(x,v)\in \set{(x,v)\in N\cC:|v|<\varepsilon|x|}\mapsto x+v\in\bb R^N,$ which is well-defined for small $\varepsilon>0$ since $\Sigma$ is compact and smooth.
	Let $\pi$ be the projection from this tubular neighborhood back to $\cC.$
	By \eqref{eq:Su-AC} with $k=0,$ the map $h:=\pi\circ F\circ\varphi$	is well-defined on $(R',\infty)\times\Sigma$ for $R'$ large enough.
	This $h$ is a diffeomorphism if $R'$ is taken even larger and leads to a graph function $u(x):= f(h^{-1}(x))-x$ that satisfies \eqref{eq:AC-rate}.
\end{proof}

Thus, we can freely use Definition~\ref{def:AC} or the intrinsic one in \cite{S20}*{Definition~2.2} when talking about asymptotically conical submanifolds.

From now on, we restrict to the Lagrangian case.
As indicated in the introduction, when we say a Lagrangian is asymptotically conical, we always assume that the asymptotic cone is a special Lagrangian cone, which consists of special Lagragians with the same angle (with no multiples of~$2\pi$ allowed).
As we always assume a positive lower bound of $\Lambda,$ by \cite{LW26}*{Lemma~6.1}, this special Lagrangian cone is in fact the union of finitely many almost isoclinic planes with the same Lagrangian angles.

\subsection{Evolution equation}

In this section, we will derive the evolution equation for $\log\Lambda$ along an almost isoclinic mean curvature flow and use it to prove the preservation of almost isoclinicity along a Lagrangian mean curvature flow.
First, we recall the following consequence of the Hessian calculations in \cite{LW26}*{Section~3.2}.

\begin{prop}
\label{prop:evol-basic}
    Let $n\ge 2$ and $M_t$ be a Lagrangian mean curvature flow of almost isoclinic submanifolds in $\bb R^{2n}$ for $t\in I\sbst\bb R.$
    Then the function $\Lambda$ satisfies
    \begin{align*}
    (\bd_t-\D)\log\Lambda
    \ge \frac 23|A|^2+
    \frac1{n^2}|\nabla\log\Lambda|^2
    \end{align*}
    for all $t\in I$ and $x\in M_t.$
\end{prop}

\begin{proof}
    Recall that the Gauss map evolves by the harmonic map heat flow along a mean curvature flow \cite{W03}. 
    Combining this with
    \cite[Proposition 3.2 and Remark 3.3]{LW26} yields
    \begin{equation}
    \label{eq:full-lambda-evolution}
    \begin{split}
    (\bd_t-\D)\log\Lambda={}&
    \sum_{\ell}\sum_{j<k}
    \pr{
    	\frac{(1+\lambda_j^2)(1+\lambda_k^2)}{(1+\lambda_j\lambda_k)^2}
    	(h_{jj\ell}-h_{kk\ell})^2+4h_{jk\ell}^2
    }
    \\[-2pt]
    &+\sum_{\ell}\sum_{j<k}\sum_{m\ne j,k}
    \frac{2(1+\lambda_j\lambda_k)(1+\lambda_m^2)}
    {(1+\lambda_j\lambda_m)(1+\lambda_k\lambda_m)}h_{jk\ell}^2.
    \end{split}
    \end{equation}
    If we let $e_\ell$'s form an local orthonormal frame near a point $p\in M,$ we can use the derivations of $dF^\Omega_S(A)$ and $dF^\Lambda_S(A)$ in the proof of \cite{LW26}*{Lemma~5.7} to derive
    \begin{align*}
	\n_{e_\ell}\pr{\log\Lambda}(p)
	= d\pr{\log\Lambda}_p(e_\ell)
	=& \sum_{j<k}
	\frac{1}{1+\lambda_j\lambda_k}
	\pr{h_{jj\ell}\lambda_k \pr{1+\lambda_j^2}
		+ h_{kk\ell}\lambda_j \pr{1+\lambda_k^2}}\\
	&- (n-1)\sum_m \frac{\lambda_m}{1+\lambda_m^2}
	\cdot h_{mm\ell}\pr{1+\lambda_m^2}\\
	=& \sum_m h_{mm\ell}
	\sum_{k\neq m}
	\pr{\frac{\lambda_k(1+\lambda_m^2)}{1+\lambda_m\lambda_k} - \lambda_m}\\
	=& \sum_m h_{mm\ell}
	\sum_{k\neq m} \frac{\lambda_k-\lambda_m}{1+\lambda_k\lambda_m}.
    \end{align*}
    This implies
	\begin{align*}
	|\n_{e_\ell}\log\Lambda|^2
	= \abs{\sum_{m,k} 
		h_{mm\ell} \frac{\lambda_k-\lambda_m}{1+\lambda_k\lambda_m}}^2
	= \abs{\sum_{j<k} 
		(h_{jj\ell}-h_{kk\ell}) \frac{\lambda_k-\lambda_j}{1+\lambda_k\lambda_j}}^2
	& \le n^2 \sum_{j<k} 
	(h_{jj\ell}-h_{kk\ell})^2 \frac{(\lambda_k-\lambda_j)^2}{(1+\lambda_k\lambda_j)^2}
	\end{align*}
    by the Cauchy--Schwarz inequality.
    Note that
    \begin{align}\label{eq:evol-first-term-CS}
    \frac{(1+\lambda_j^2)(1+\lambda_k^2)}{(1+\lambda_j\lambda_k)^2}
    = 1 + \frac{(\lambda_k-\lambda_j)^2}{(1+\lambda_k\lambda_j)^2},
    \end{align}
    so among the terms on the right hand side of \eqref{eq:full-lambda-evolution}, it remains to estimate
    \[
    Q_0=\sum_{\ell,j<k}(h_{jj\ell}-h_{kk\ell})^2
    +4\sum_{\ell,j<k}h_{jk\ell}^2
    \] 
    by $|A|^2.$
    The full symmetry of $h_{jk\ell}$ gives the exact grouping
    \begin{align*}
    	Q_0={}&\sum_\ell
    	\sum_{j\ne\ell}\pr{(h_{\ell\ell\ell}-h_{jj\ell})^2+4h_{jj\ell}^2}
    	+ \sum_\ell\sum_{\substack{j<k\\j,k\ne\ell}}(h_{jj\ell}-h_{kk\ell})^2 
    	+12\sum_{j<k<\ell}h_{jk\ell}^2.
    \end{align*}
    Since
    \(
    (a-b)^2+4b^2-\tfrac23(a^2+3b^2)=\tfrac13(a-3b)^2\ge0
    \)
    and
    \begin{align*}
	|A|^2
	= \sum_j h_{jjj}^2
	+ \sum_{j\neq k}3h_{jjk}^2 
	+ \sum_{j<k<\ell}6h_{jk\ell}^2,
	\end{align*}
    we obtain $Q_0\ge\tfrac23|A|^2$ using $n-1\ge1$. 
    This proves the proposition. 
\end{proof}

\begin{cor}
\label{cor:preserve-AI}
    Let $n\ge 3$ and $M$ be a complete, asymptotically conical, zero-Maslov, almost isoclinic Lagrangian with $\Lambda\ge c>0.$
    Then there exists a maximal existence time $T\in (0,\infty]$ such that there exists a unique smooth Lagrangian mean curvature flow $M_t$ for $t\in[0,T)$ with $M_0=M$ satisfying that each $M_t$ is complete, asymptotically conical, and almost isoclinic with $\Lambda\ge c.$
\end{cor}

Here we assume $n\ge 3$ since we will apply the short-time existence result in \cite{S20}.
Recall that when we say a Lagrangian is asymptotically conical, we always assume that the asymptotic cone is a special Lagrangian cone.

\begin{proof}
	The uniqueness part follows from \cite{CY07}. 
	For the existence, first, we apply \cite{S20}*{Theorem~4.3} to get a smooth asymptotically conical $M_t$ for $t\in[0,T)$ where $T\in(0,\infty]$ is the maximal existence time of the flow and $M_0=M.$
    By \cite{S20}*{Proposition 4.6}, for each $t\in[0,T),$ $M_t$ is asymptotically conical with the same rate and the same cone as $M_0.$
    
    It remains to show that each $M_t$ is almost isoclinic with $\Lambda\ge c.$
    Fix any compact interval $[0,t_0]\sbst[0,T).$
    The asymptotic conicality condition allows us to write $M_t\setminus \ovl B_R$ as the graph of a function with a decay rate $\beta-1-k<-1-k$ for all $k$-th derivatives according to Definition~\ref{def:AC}.
    Thus, the bounds on the initial derivatives give uniform continuity in time of the Hessian for a short interval.
    Since the ordered eigenvalues of a real symmetric matrix depend continuously on its entries, the uniform continuity of the Hessian of the graph function in time and the uniform bound $\Lambda\ge c$ imply that
    $1+\lambda_j\lambda_k>c/2$ for all $j\ne k$ in an almost isoclinic chart on a sufficiently short time interval $[0,\varepsilon]\sbst[0,t_0].$
    The almost isoclinicity is preserved on $M_t\cap \ovl B_R$ on $[0,\varepsilon]$ after possibly shrinking $\varepsilon$ by the compactness.
    Thus, it follows that $M_t$ is almost isoclinic for $t\in[0,\varepsilon],$ and hence Proposition~\ref{prop:evol-basic} implies that $\Lambda\ge c$ on $M_t$ for $t\in[0,\varepsilon].$
    This $\varepsilon$ depends only on the asymptotic bounds and the bound $c,$ which are both preserved, so both the almost isoclinicity and the property $\Lambda\ge c$ are preserved for all $t\in[0,T).$    
\end{proof}

Another consequence of the monotonicity formula given by Proposition~\ref{prop:evol-basic} is a classification results for tangent flows.
The argument is adapted from the main theorem in \cite{N07} (see also \cite{W01-JDG}*{Remark 5.1}) based on the better reaction term given by the evolution of $\log\Lambda.$

\begin{cor}
\label{cor:tangent-flow-structure}
	Let $M_t$ be an almost isoclinic Lagrangian mean curvature flow in $\bb R^{2n}$ for $t\in[0,T)$ with $\Lambda\ge c>0$ and finite entropy.
    Suppose there is a globally defined Lagrangian angle on $M_t$ which is uniformly bounded.
	Given a singularity $(x_0, T),$ a tangent flow of $M_t$ at $(x_0, T)$ is a finite union of special Lagrangian planes.
	
	To be more precise, given a sequence $\lambda_j\to \infty,$ consider the sequence of rescaled flows
	\begin{align*}
		M^j_t
		:= \lambda_j\pr{M_{T+t\lambda_j^{-2}} - x_0}
	\end{align*}
	for $t\in [-\lambda_j^2T,0).$
	Then there exist numbers $\ovl \theta_1,\cdots,\ovl\theta_N$ and special Lagrangian planes $L_1,\cdots,L_N$ such that the following holds.
	After passing to a subsequence of $\lambda_i,$ we have that for any compactly supported $\varphi\in C^\infty\pr{\bb R^{2n}},$ any $f\in C^2(\bb R),$ and any $t<0,$
	\begin{align*}
		\lim_{j\to\infty}
		\int_{M^j_t} f\pr{e^{i\theta^j_t}}\varphi\,dx
		= \sum_{k=1}^N m_k\cdot f\pr{e^{i\ovl\theta_k}}\cdot \mu_k(\varphi)
	\end{align*}
	where $\mu_k$'s and $m_k$'s are the Radon measures of the support of $L_k$'s and the multiplicities of $L_k$'s.
\end{cor}

\begin{proof}
	By Proposition~\ref{prop:evol-basic},
	$(\bd_t-\D)\log \Lambda\ge (2/3)|A|^2.$
	Thus, if we let
	\begin{align}\label{eq:b-heat}
	\Phi_{x_0,T}(x,t)
	:= \frac 1{4\pi(T-t)^{n/2}} 
	e^{-|x-x_0|^2/4(T-t)},
	\end{align}
	by Huisken's monotonicity formula, it follows that
	\begin{align*}
		\bd_t\int_{M_t} \pr{-\log\Lambda}
		\cdot \Phi_{x_0,T}\,dx
		& = \int_{M_t}
		\pr{(\bd_t-\D) \pr{-\log\Lambda}
			- \abs{\h + \frac{(x-x_0)^\perp}{2(T-t)}}^2 \pr{\log\Lambda}^2}\Phi_{x_0,T}\,dx\\
		& \le -\int_{M_t} 
		\pr{\frac 23|A|^2 + \abs{\h + \frac{(x-x_0)^\perp}{2(T-t)}}^2 \pr{\log\Lambda}^2}\Phi_{x_0,T}\, dx
	\end{align*}
	and
	\begin{align*}
		\bd_t\int_{M_t} \Phi_{x_0,T}\,dx
		& = -\int_{M_t}
		\abs{\h + \frac{(x-x_0)^\perp}{2(T-t)}}^2 \Phi_{x_0,T}\,dx.
	\end{align*}
	These monotonicity formulas then imply
	\begin{align*}
		\lim_{t\to T}
		\int_t^T \int_{M_s} 
		\pr{\frac 23|A|^2 + \abs{\h + \frac{(x-x_0)^\perp}{2(T-t)}}^2}\Phi_{x_0,T}\,dx\,ds=0.
	\end{align*}
	In particular, 
	\begin{align*}
		\lim_{t\to T}
		\int_t^T \int_{M_s}
		|\h|^2 \Phi_{x_0,T}\,dx\,ds=0,
	\end{align*}
	so we also derive
	\begin{align*}
		\lim_{t\to T}
		\int_t^T \int_{M_s}
		\abs{\frac{(x-x_0)^\perp}{2(T-t)}}^2 \Phi_{x_0,T}
		\,dx\,ds=0.
	\end{align*}
	
	For the rescaled flows, it follows that for any $a<b<0,$ 
	\begin{align}\label{lim-at-time-t}
		\lim_{j\to\infty} \int_a^b \int_{M^i_t} \pr{|A|^2 + |x^\perp|^2} \Phi \,dx\,dt=0.
	\end{align}
	We can then choose $a<0$ such that
	\begin{align}\label{limit-flat}
		\lim_{j\to\infty} \int_{M^i_a\cap B_R} \pr{|A|^2 + |x^\perp|^2} \Phi \,dx=0
	\end{align}
	for all $R>0.$
	Since the initial submanifold has finite entropy, there exists a constant $C<\infty$ such that
	\begin{align*}
		\cH^n\pr{M^j_a\cap B_R}
		\le CR^n
	\end{align*}
	for all $j\in\bb N$ and $R>0.$
	We can then apply the compactness theorem \cite[Proposition~5.1]{N07} to $M^j_a$ such that after passing to a subsequence, there exist integral special Lagrangian cones $L_1,\cdots,L_n$ and numbers $\ovl\theta_1,\cdots,\ovl\theta_N$ such that for any smooth function $\varphi$ with compact support and any $f\in C^2(S^1),$
	\begin{align}\label{limit-at-time-a}
		\lim_{j\to\infty}
		\int_{M^j_a}f\pr{e^{i\theta^j_a}}\varphi\,dx
		= \sum_{k=1}^n
		m_k\cdot f\pr{e^{i\ovl\theta_k}}\mu_k\pr{\varphi}
	\end{align}
	where $\mu_k$ and $m_k$ are the Radon measure and the multiplicity of $L_k.$

	We now deduce from \eqref{limit-flat} that each $L_k$ is a union of special Lagrangian planes with angle $\ovl\theta_k.$
	By \cite{H86}*{Theorem 5.3.2}, the limit of $M^j_a$ is a curvature varifold $V$ and satisfies
	\begin{align*}
	\int_{B_R}|A_V|^2 d\|V\|
	\le \liminf_{j\to\infty} \int_{M^i_a\cap B_R} |A|^2\, dx
	=0
	\end{align*}
	for all $R>0,$ where $A_V$ is the generalized second fundamental form 
    \footnote{In \cite{H86}, the generalized second fundamental form of a varifold $V$ in an open set $E$ is denoted by $\mathscr S_EV;$ in~\cite{A26}, the generalized second fundamental form of a varifold $V$ is denoted by $\mathcal C_V.$ }
    of $V$ and the equality is by~\eqref{limit-flat}.
	This implies $A_V$ is $\|V\|$-almost everywhere zero, and hence \cite{A26}*{Theorem~4.3} implies that $V$ is locally given by the sum of multiplicity planar measures near the origin.
	The conicality then implies that $V$ is the sum of multiplicity planar measures, and we put those with the same Lagrangian angles together and form $L_k$'s, which are now known to be a special Lagrangian cone which is a union of planes.

	Next, we will show that we can extend \eqref{limit-at-time-a} to any negative time $b<0.$  
	From \cite{S02}, for any $j$ and $t,$
	$(\bd_t-\D) e^{i\theta^j_t}
	= e^{i\theta^j_t} |\h|^2.$
	Thus, for a given $f\in C^2(S^1),$ if we write $f=f\pr{e^{i\theta^j_t}},$ then
	\begin{align*}
		\abs{(\bd_t-\D) f}
		+ \abs{\n f}^2
		& = \abs{f'
			\cdot (\bd_t-\D) e^{i\theta^j_t}
			-f'' \abs{\n e^{i\theta^j_t}}^2}
		+ \abs{f'\cdot \n e^{i\theta^j_t}}^2\\
		& = \abs{f'\cdot e^{i\theta^j_t}|\h|^2
			- f'' \cdot |\h|^2}
		+ |f'|^2|\h|^2\\
		& \le C_f \cdot |\h|^2
	\end{align*}
	for some $C_f<\infty$ depending on $f.$
	Thus,
	\begin{align*}
		\bd_t \int_{M^j_t}
		f\cdot \varphi\,dx
		& = \int_{M^j_t}
		\pr{
			(\bd_t-\D)f\cdot \varphi
			+ f\pair{\h, D\varphi}
			- \pair{\n f, \n\varphi}
			- f|\h|^2\varphi
		}\,dx
	\end{align*}
	and hence for any $b<0,$
	\begin{align*}
		\abs{\int_{M^j_b}
			f\cdot \varphi\,dx
			- \int_{M^j_a}
			f\cdot \varphi\,dx}
		& = \abs{\int_a^b \int_{M^j_t}
			\pr{
				(\bd_t-\D)f\cdot \varphi
				+ f\pair{\h, D\varphi}
				- \pair{\n f, \n\varphi}
				- f|\h|^2\varphi
			}\,dx\,dt}\\
		& \le \int_a^b \int_{M^j_t}
		\pr{C_f|\h|^2 \varphi
			+ |f|\cdot |\h|\cdot |D\varphi|
			+ \sqrt{C_f}|\h|\cdot |\n \varphi|
			+ |f\varphi|\cdot |\h|^2
		}\,dx\,dt\\
		& \le C \pr{
			\int_a^b \int_{M^j_t} |\h|^2\Phi\,dx\,dt
			+ \pr{\int_a^b \int_{M^j_t} |\h|^2\Phi\,dx\,dt}^{1/2}
		}
	\end{align*}
	for some $C<\infty$ depending on $f,\varphi,a,$ and $b.$
	Thus, by \eqref{lim-at-time-t}, we get
	\begin{align*}
		\lim_{j\to\infty} \int_{M^j_b}
		f\cdot \varphi\,dx
		= \lim_{j\to\infty} \int_{M^j_a}
		f\cdot \varphi\,dx
		= \sum_{k=1}^n
		m_k\cdot f\pr{e^{i\ovl\theta_k}}\mu_k\pr{\varphi}.
	\end{align*}
	This finishes the proof.	
\end{proof}

As a corollary of this classification of tangent flows, we can rule out type-I singularities in finite time along such a Lagrangian mean curvature flow.
Such a result has been established in the almost calibrated case \cite{W01-JDG}.

\begin{cor}
    Let $M_t$ be a zero-Maslov, asymptotically conical, almost isoclinic Lagrangian mean curvature flow in $\bb R^{2n}$ with $\Lambda\ge c>0$ and finite entropy.
    Then $M_t$ does not admit a finite-time type-I singularity.
\end{cor}

Another immediate corollary is the long-time existence of Lagrangian mean curvature flows starting from low-entropy submanifolds.

\begin{cor}
    Let $M_t$ be a zero-Maslov, asymptotically conical, almost isoclinic Lagrangian mean curvature flow in $\bb R^{2n}$ with $\Lambda\ge c>0.$
    If $\lambda(M_0)\le 2,$ then $M_t$ exists for all $t\ge 0.$
\end{cor}

\section{\bf Graphical case}

\subsection{Graphicality criterion and consequence}
\label{sec:G-criterion}

In this section, we generalize \cite{LW26}*{Proposition~6.3} and get a graphicality criterion for almost isoclinic Lagrangians which is ``calibrated enough.''

\begin{prop}
\label{prop:graphicality}
    Given any $c>0,$ there exists $\delta>0$ such that the following holds.
    Suppose $M$ is a complete, connected, almost isoclinic Lagrangian in $\bb R^{2n}$ with a globally defined Lagrangian angle $\theta.$
	If $\Lambda\ge c$ and if there exists $\theta_0\in\bb R$ such that
    \begin{align*}
    |\theta-\theta_0|<\delta
    \end{align*}
    on $M,$ then $M$ is graphical up to a rotation.
    To be precise, there exist $\varphi\in\bb R$ and a (vector-valued) smooth function $V\colon \bb R^n\to \bb R^n$ with $|\n V|\le C<\infty$ such that
	\begin{align*}
	e^{i\varphi} M = {\rm Graph}_{\bb R^n} V.
	\end{align*}
\end{prop}

\begin{proof}
    For any $p\in M,$ there exist $\varphi\in\bb R$ and an $n$-plane $L_\varphi=e^{i\varphi}L_B$ so that $T_pM$ is graphical over $L_\varphi$ and satisfies the condition in Definition~\ref{def:AI}.
    Write $T_pM$ as the graph of a symmetric matrix $S_\varphi$ over $L_0,$ and let $\lambda_k$'s be the eigenvalues of $S_\varphi.$
    Thus, if we let $\theta_k:=\arctan \lambda_k$ be the characteristic angles in this chart, it follows that
    \begin{align}\label{SL-ang}
    e^{i\theta(p)} = e^{i\pr{n\varphi + \sum_k\theta_k}}.
    \end{align}
    Consider
    \begin{align*}
    \psi:=\varphi + \frac 1n \sum_k\theta_k,
    \end{align*}
    and hence \eqref{SL-ang} becomes $e^{i\theta(p)} = e^{in\psi}.$
    Moreover, if we look at the new reference plane $L_\psi=e^{i\psi}L_B,$ according to \cite{LW26}*{Lemma~2.6}, $T_pM$ is graphical over $L_\psi$ and the eigenvalues of the graph matrix are
    \begin{align*}
    \td \lambda_k 
    = \frac{\lambda_k\cos\ovl\theta - \sin\ovl\theta}{\cos\ovl\theta + \lambda_k\sin\ovl\theta}
    = \tan(\theta_k-\ovl\theta)
    \end{align*}
    where we write $\ovl\theta:=\frac 1n \sum_k\theta_k.$
    Note that each $\theta_k-\ovl\theta$ still lies in $(-\pi/2,\pi/2),$ so the new characteristic angles in this chart are
    $\td\theta_k:=\arctan\td\lambda_k=\theta_k-\ovl\theta.$

    Now, we take a positive number $\delta$ such that
    \begin{align*}
    \delta<\frac\pi 2 - \frac{n-1}n\arccos c.
    \end{align*}
    We claim that $T_pM$ is also graphical over $e^{i(\psi+(\theta_0-\theta(p))/n)}L_B$ if $|\theta(p)-\theta_0|<\delta.$
    Since the angles satisfy
    \begin{align*}
    \Lambda = \prod_{j<k}\cos\pr{\td\theta_j-\td\theta_k}\ge c,
    \end{align*}
    each $\cos\pr{\td\theta_j-\td\theta_k}$ is bounded by $c,$ so 
    \begin{align*}
    |\td\theta_j-\td\theta_k|\le \arccos c\in \left[0, \frac\pi 2\right).
    \end{align*}
    Then the relation $\sum_k\td\theta_k=0$ implies 
    \begin{align}\label{theta-k-bdd}
    |\td\theta_k|\le \frac{n-1}n\arccos c.
    \end{align}
    As a result, by \cite{LW26}*{Lemma~2.6}, a characteristic angle in the chart  $e^{i(\psi+(\theta_0-\theta(p))/n)}L_B$ is given by the arctangent of
    \begin{align}\label{eq:new-angle}
    \lambda'_k
    = \frac{\td\lambda_k\cos(\theta_0-\theta(p)) - \sin(\theta_0-\theta(p))}{\cos(\theta_0-\theta(p)) + \td\lambda_k\sin(\theta_0-\theta(p))}
    = \tan\pr{\td\theta_k-(\theta_0-\theta(p))},
    \end{align}
    which is still finite by $|\theta(p)-\theta_0|<\delta,$ the definition of $\delta,$ and the bound \eqref{theta-k-bdd}.
    Since
    \begin{align*}
    e^{i\pr{n\psi + \pr{\theta_0-\theta(p)}}}
    = e^{i\theta_0},
    \end{align*}
    which is independent of $p,$ it follows that there is an open cover of $M$ given by
    \begin{align*}
    M\sbst \bigcup_{t=0}^{n-1}
    \mathcal O^M_{\pr{\theta_0+2\pi t}/n}
    \end{align*}
    where we let $\mathcal O^M_\alpha
    := \set{q\in M:T_qM\in \mathcal O_\alpha}.$

    We claim that each piece $\mathcal O^M_{\pr{\theta_0+2\pi t}/n}$ is a closed subset in $M.$
    Fix a $t\in\set{0,1,\cdots,n-1}.$
    Suppose $p\in \mathcal O^M_{\pr{\theta_0+2\pi t}/n}$ and let
    $\td\theta_k -(\theta_0-\theta(p))$
    be the characteristic angles in this chart as derived in \eqref{eq:new-angle}.
    By $|\theta(p)-\theta_0|<\delta,$ the definition of $\delta,$ and the bound \eqref{theta-k-bdd}, such an angle satisfies
    \begin{align}\label{eq:new-grad-bdd}
    |\td\theta_k -(\theta_0-\theta(p))|
    \le |\td\theta_k|+\delta
    \le \frac{n-1}n\arccos c+\delta,
    \end{align}
    which is a constant strictly less than $\pi/2.$    
    This implies that if $p_\ell$ is a sequence of points in $\mathcal O^M_{\pr{\theta_0+2\pi t}/n}$ that converges to some $p_\infty\in M,$ then the characteristic angles of $T_{p_\infty}M$ still satisfy the same estimates \eqref{eq:new-grad-bdd}, so $p_\infty\in O^M_{\pr{\theta_0+2\pi t}/n}.$
    This proves that each $\mathcal O^M_{\pr{\theta_0+2\pi t}/n}$ is closed.
    The connectivity of $M$ then implies 
    \begin{align*}
    M = \mathcal O^M_{\pr{\theta_0+2\pi t_0}/n}
    \end{align*}
    for a common $t_0\in\set{0,1,\cdots,n-1}.$ 
    This and the gradient bound \eqref{theta-k-bdd} allow us to conclude that the orthogonal projection 
    \begin{align*}
    \pi\colon \bb R^{2n}\to L_M:= e^{i\pr{\theta_0+2\pi t_0}/n} L_B
    \end{align*}
    restricts to a local diffeomorphism on $M.$

    The rest of the proof is the same as \cite{LW26}*{Proposition~6.3}.
    By the gradient bound \eqref{eq:new-grad-bdd} and a covering map argument, it follows that $\pi|_M$ is a global diffeomorphism and the conclusion follows.
\end{proof}

As a consequence, we show that an immortal almost isoclinic Lagrangian mean curvature flow will eventually become graphical if it admits a globally defined Lagrangian angle.
This result will allow us to reduce all the convergence results to the graphical case once we establish their long-time existence.

\begin{cor}
\label{cor:long-time-existence-imply-conv}
    Suppose $M$ is a complete, connected, almost isoclinic Lagrangian in $\bb R^{2n}$ with a globally defined Lagrangian angle and suppose $M$ is asymptotic to a union of almost isoclinic Lagrangian $n$-planes with the same Lagrangian angles.
    If the mean curvature flow $M_t$ starting from $M$ is smoothly defined for all $t\ge 0,$ then there exists $T<\infty$ and $\varphi\in\bb R$ such that for all $t\ge T,$ $e^{i\varphi}M_t$ is a graphical mean curvature flow over the fixed Lagrangian base plane. 
\end{cor}

In particular, as a consequence of this corollary, each $M_t$ along the flow only has exactly one end and is diffeomorphic to $\bb R^n.$

\begin{proof}
	By the condition on the asymptotic planes, there exists $c>0$ such that $\Lambda\ge c$ on $M_t$ for all $t\ge 0.$	
	Let $\delta$ be the constant from Proposition~\ref{prop:graphicality} that depends on $c.$
	We will use the heat kernel
	\begin{align*}
	\Phi(x,t):= \frac 1{t^{n/2}} e^{-|x|^2/4t}
	\end{align*}
	to show that the variation of the Lagrangian angle $\theta$ of $M_t$ will decrease to a desired range.
	Note that the heat kernel satisfies 
	\begin{align}\label{eq:heat-kernel}
	(\bd_t-\D)\Phi = \frac{|x^\perp|^2}{4t^2}\Phi
	\end{align}
	on $M_t$ for $t>0.$
	The commonly stated version is for the backward heat kernel (e.g., \eqref{eq:b-heat}), and for \eqref{eq:heat-kernel}, it follows from
	\begin{align*}
	\bd_t \Phi - \D\Phi 
	= \pr{-\frac n{2t} + \frac{|x|^2}{4t} - \frac{\pair{x,\h}}{2t}}\Phi
	- \pr{-\frac n{2t} + \frac{|x^T|^2}{4t} - \frac{\pair{x,\h}}{2t}}\Phi
	= \frac{|x^\perp|^2}{4t^2}\Phi.
	\end{align*}
	
	Let $\theta_0$ be the Lagrangian angle of those asymptotic planes.
	At time $t=1,$ $M_1$ is asymptotically flat, so there exists $R>0$ such that $|\theta(\cdot,1)-\theta_0|\le \delta/2$ on $M_1\setminus B_R.$
	Take $C>0$ such that
	\begin{align*}
	|\theta(\cdot,1)-\theta_0|\le \frac\delta 2 + C\cdot \Phi(\cdot,1)
	\end{align*}
	on the whole $M_1.$
	Combining $(\bd_t-\D)\theta=0$ from \cite{S02} with \eqref{eq:heat-kernel} implies
	\begin{align*}
	\pr{\bd_t-\D}
	\pr{\theta-\frac\delta 2 - \theta_0-C\cdot\Phi}
	= -\frac{C|x^\perp|^2}{4t^2}\Phi
	\le 0.
	\end{align*}
	The maximum principle with the asymptotic flatness implies 
	\begin{align*}
	\theta-\theta_0\le \frac\delta 2 + C\cdot \Phi
	\end{align*}
	on $M_t$ for all $t\ge 1.$
	The same argument implies
	\begin{align*}
	-\pr{\theta-\theta_0}\le \frac\delta 2 + C\cdot \Phi
	\end{align*}
	on $M_t$ for all $t\ge 1.$ 
	Combining them leads to
	\begin{align*}
	|\theta-\theta_0|
	\le \frac\delta 2 + C\cdot\Phi
	\le \delta
	\end{align*}
	on $M_t$ for $t\ge T_\delta$ where $T_\delta\in(0,\infty)$ can be taken to be $(2C/\delta)^{-n/2}.$
	Proposition~\ref{prop:graphicality} then implies that $M_t$ is graphical for $t\ge T_\delta.$
\end{proof}

\subsection{Graphical mean curvature flow}
\label{sec:G-MCF}

In this section, we establish almost isoclinic analogues of the three main results of Tsai--Tsui--Wang \cite{TTW25}. 
We prove long-time existence for entire gradient graphs with bounded initial Hessian, rescaled convergence in the bounded-gradient case, and convergence to the self-expander determined by an initial blow-down. 
The proofs use the evolution equations for two geometric quantities $\omega$ and $\Lambda$, see \eqref{eq:w and Lambda}, \eqref{eq:two functions}, a localized maximum principle, and the approximation lemma in Appendix~\ref{sec:approximation}.

Throughout, $n\ge2$, and all pairwise eigenvalue conditions refer to \emph{distinct} indices. 
For a real symmetric matrix $B$, set
\[
 \Theta(B)=\sum_{j=1}^n\arctan\lambda_j(B).
\]
We use this real-valued branch of the Lagrangian angle. 
The potential equation \footnote{This is equivalent to equation (1.1) in \cite{TTW25} } is
\begin{equation}\label{eq:flow}
 u_t=\Theta(D^2u),\qquad u(x,0)=u_0(x),
 \qquad (x,t)\in\R^n\times[0,\infty).
\end{equation}
The graph $L_t=\{(x,Du(x,t)):x\in\R^n\}\subset\C^n$ evolves by mean curvature up to a tangential reparametrization.

A potential is \emph{weakly almost isoclinic} in this graph chart when $1+\lambda_j\lambda_k\ge0$ for every $j\ne k$, and is \emph{strictly almost isoclinic} when all these inequalities are strict. 
Two-convexity (see Definition 2.1 in \cite{TTW25}) imposes the additional condition $\lambda_j+\lambda_k\ge0$.
Write
\begin{equation}\label{eq:w and Lambda}
 \omega(B)=\prod_{j=1}^n\frac{1}{\sqrt{1+\lambda_j^2}}
 \,\,\text{ and }\,\,
 \Lambda(B)=\prod_{j<k}
 \frac{1+\lambda_j\lambda_k}{\sqrt{(1+\lambda_j^2)(1+\lambda_k^2)}}.
\end{equation}
On the strict almost isoclinic region, $0<\Lambda\le1$. 
We write
\begin{equation}\label{eq:two functions}
\omega[u]=\omega(D^2u)
\,\,\text{ and }\,\,
\Lambda[u]=\Lambda(D^2u).\end{equation}

The main results of this section are the following three theorems.
The first one is a long-time existence result.

\begin{thm}
\label{thm:existence}
Let $u_0\in C^2(\R^n)$ satisfy
\[
 \sup_{\R^n}|D^2u_0|^2\le c
\,\,\text{ and }\,\,
 1+\lambda_j(D^2u_0)\lambda_k(D^2u_0) \ge0.
\]
Here and below, $|D^2u|$ denotes the Euclidean Frobenius norm. 
Then \eqref{eq:flow} has a solution
\[
 u\in C^0(\R^n\times[0,\infty))
       \cap C^\infty(\R^n\times(0,\infty)),
\]
unique among such solutions whose Hessians are uniformly bounded on finite time intervals. 
It has the following properties.
\begin{enumerate}[label=\textup{(\roman*)}]
\item Weak almost isoclinicity is preserved. 
If $u_0$ is pointwise strictly almost isoclinic, then $u(\cdot,t)$ is strictly almost isoclinic for every $t>0$, without assuming a uniform positive initial gap.
\item If $\omega_0=\inf_{\R^n}\omega[u_0]>0$, then $\omega[u](x,t)\ge\omega_0$. 
In particular, the Hessian remains uniformly bounded.
\item For every integer $\ell\ge2$ there is $C_\ell=C_\ell(n,c)<\infty$ such that
\begin{equation}\label{eq:derivative-decay}
 \sup_{\R^n}|D^\ell u(\cdot,t)|^2\le C_\ell t^{2-\ell}
\end{equation}
for all $t>0.$
Moreover,
\begin{equation}\label{eq:time-displacement}
 |u(x,t)-u_0(x)|\le\frac{n\pi}{2}t.
\end{equation}
\end{enumerate}
\end{thm}

The next theorem asserts a rescaled limit under a bounded gradient assumption.

\begin{thm}
\label{thm:rescaled}
Under the hypotheses of Theorem~\ref{thm:existence}, suppose additionally $\sup_{\R^n}|Du_0|\le L<\infty$. 
Then $\sup_{\R^n}|Du(\cdot,t)|\le L$ for all $t>0.$
The entire rescaled family
\begin{equation}\label{eq:potential-rescaling}
 u_\mu(x,t)=\mu^{-2}u(\mu x,\mu^2t)
\end{equation}
converges to zero in $C^\infty_{\mathrm{loc}}(\R^n\times(0,\infty))$ as $\mu\to\infty$. 
Equivalently, for the rescaled time variable $s=\log\sqrt{2t+1},$
the normalized Lagrangian graphs
\[
 \widetilde L_s:=\frac 1{\sqrt{2t+1}} L_t
 =\left\{\left(y,\frac{Du(\sqrt{2t+1}\, y,t)}{\sqrt{2t+1}}\right):y\in\R^n\right\}
\]
converge locally smoothly to $\R^n\times\{0\}$. 
In fact, their height satisfies
\begin{equation}\label{eq:normalized-height}
 \sup_y\frac{|Du(\sqrt{2t+1}\, y,t)|}{\sqrt{2t+1}}\le\frac{L}{\sqrt{2t+1}}.
\end{equation}
\end{thm}

The following theorem assumes the existence of a full blow-down limit in \eqref{eq:initial-blowdown} and asserts a convergence to an expander.

\begin{thm}
\label{thm:expander}
Under the hypotheses of Theorem~\ref{thm:existence}, assume that for every $x\in\R^n,$ the limit
\begin{equation}\label{eq:initial-blowdown}
 \lim_{\mu\to\infty}\mu^{-2}u_0(\mu x)=U_0(x)
\end{equation}
exists. 
Then $U_0\in C^{1,1}(\R^n)$, the convergence is in $C^1_{\mathrm{loc}}$, and $U_0$ is homogeneous of degree two. 
The whole family \eqref{eq:potential-rescaling} converges in $C^\infty_{\mathrm{loc}}(\R^n\times(0,\infty))$ to a solution $U$ of \eqref{eq:flow} satisfying $U(\cdot,t)\to U_0$ locally uniformly as $t\searrow 0$. 
It is weakly almost isoclinic for positive time and satisfies
\begin{equation}\label{eq:selfsimilar}
 U(\rho x,\rho^2t)=\rho^2U(x,t),\,
 U(x,t)=tV(x/\sqrt t),\,\text{ and }\, 
 V=U(\cdot,1).
\end{equation}
Thus,
\begin{equation}\label{eq:expander-equation}
 \Theta(D^2V)=V-\tfrac12x\cdot DV,
\end{equation}
and the Lagrangian graphs are self-expanding. 
\end{thm}

\subsubsection{Evolution and localization}
The induced metric of the gradient graph is $g=I+(D^2u)^2$. In a frame
diagonalizing $D^2u$, put $\theta_j=\arctan\lambda_j.$
Write
\[
 e_j=\frac{a_j+\lambda_jJa_j}{\sqrt{1+\lambda_j^2}}
 \,\,\text{ and }\,\,
 h_{jk\ell}=\ip{A(e_j,e_k)}{Je_\ell}.
\]
The cubic tensor $h$ is completely symmetric, and $|A|^2=\sum_{j,k,\ell}h_{jk\ell}^2$.

\begin{lemma}
\label{lem:evolution}
Along a smooth strictly almost isoclinic Lagrangian graphical mean curvature flow,
\begin{align}
(\bd_t-\D) \log\omega &\ge\frac1n|\nabla\log\omega|^2, \text{ and}
 \label{eq:omega-evolution}\\
(\bd_t-\D) \log\Lambda &\ge\frac23|A|^2+
                         \frac1N|\nabla\log\Lambda|^2.
 \label{eq:lambda-evolution}
\end{align}
The first inequality remains valid under weak almost isoclinicity.
\end{lemma}

\begin{proof}
The standard evolution identity for the graph volume density gives
\begin{align*}
(\bd_t-\D)\log\omega \ge{}&\sum_j(1+\lambda_j^2)h_{jjj}^2
 +\sum_{j\ne k}(3+\lambda_j^2+2\lambda_j\lambda_k)h_{jjk}^2
+\sum_{j<k<\ell}(6+2\lambda_j\lambda_k+2\lambda_k\lambda_\ell+2\lambda_\ell\lambda_j)h_{jk\ell}^2.
\end{align*}
Also $\nabla_{e_k}\log\omega=-\sum_j\lambda_j h_{jjk}$, so $|\nabla\log\omega|^2\le n\sum_{j,k}\lambda_j^2h_{jjk}^2$.
After retaining $n^{-1}|\nabla\log\omega|^2$, all remaining coefficients
are non-negative because $\lambda_j\lambda_k\ge-1$.
This is precisely the calculation in \cite{TTW25}*{Proposition 4.1}; no
condition on $\lambda_j+\lambda_k$ is used.

The second inequality is proven in Proposition~\ref{prop:evol-basic}.
\end{proof}

We recall the following maximum principle from \cite{TTW25}*{Lemma~5.1}.

\begin{lemma}
	\label{lem:cutoff}
Let $F$ parametrize a smooth family of properly immersed submanifolds evolving by mean curvature in Euclidean space.
Suppose a positive function~$v$ satisfies
\[
 (\bd_t-\D)\log v\le-q|\nabla\log v|^2
\]
for some $q>0.$
If $v$ is initially bounded, then
$\sup_{L_t}v\le\sup_{L_0}v$.
\end{lemma}

\begin{cor}
	\label{cor:strictness}
For a strictly almost isoclinic flow, set
\[
m_R=\min_{|y|\le R}\Lambda[u_0](y).
\]
Wherever $\phi(x,t)=R^2-|x|^2-|Du(x,t)|^2-2nt>0$, one has
\begin{equation}\label{eq:local-strictness}
 \Lambda[u](x,t)\ge m_R\left(\frac{\phi(x,t)}{R^2}\right)^N.
\end{equation}
The corresponding estimate for $\omega$ has exponent $n$ and initial minimum $\min_{B_R}\omega[u_0]$.
\end{cor}

\begin{proof}
Use $v=\Lambda^{-1}$, $q=1/N$ in Lemma~\ref{lem:cutoff} and
\eqref{eq:lambda-evolution}. 
On the initial cutoff support, $|y|\le R$, so $\sup\phi_+v^q\le R^2m_R^{-1/N}$. Rearranging proves the claim.
For $\omega$, use \eqref{eq:omega-evolution} and $q=1/n$.
\end{proof}

\subsubsection{Long-time existence and derivative estimates}
We next state explicitly the analytic input used in the existence argument. 
Smooth initial potentials with all derivatives of order at least two bounded admit short-time solutions with the same property; see \cite{TTW25}*{Proposition 3.2} and \cite{CCY13}. 
The interior theorem of Nguyen--Yuan \cite{NY11}*{Theorem 1.1} says that a Hessian bound on a unit parabolic cylinder gives interior bounds for $D^3u$ and $Du_t$, provided, when $n\ge4$,
\[
 3+\lambda_j^2+2\lambda_j\lambda_k\ge0\quad\hbox{for all }j,k.
\]
The constants depend only on the dimension and the Hessian bound. 
Under weak almost isoclinicity, for $j\ne k$,
\begin{equation}\label{eq:ny-condition}
 3+\lambda_j^2+2\lambda_j\lambda_k\ge1+\lambda_j^2>0;
\end{equation}
for $i=j$ the condition is automatic. 
Thus, this theorem applies in every dimension considered here. 
Differentiation of \eqref{eq:flow} and interior Schauder estimates \cite{Friedman} then give every higher-order interior estimate. 
No positive lower bound for $1+\lambda_j\lambda_k$ enters these constants.

We start from smooth initial data.
A general short-time existence is proven in Corollary~\ref{cor:preserve-AI}.
For the graphical case, a long-time existence result follows from a Hessian bound.

\begin{prop}
	\label{prop:smooth}
Suppose $u_0$ is smooth, all $D^\ell u_0$ are bounded for $\ell\ge 2,$ $\norm{D^2u_0}_{\op}\le M,$ and for all $j\neq k,$ $1+\lambda_j\lambda_k\ge c>0.$
Then its solution exists for all positive time, preserves $\inf\omega[u_0]$ and $\inf\Lambda[u_0]$, and obeys \eqref{eq:derivative-decay}, with constants depending only on $n$ and $M$.
\end{prop}

\begin{proof}
Short-time existence and the bounds on the initial derivatives give uniform continuity in time of the Hessian for a short interval.
Since the ordered eigenvalues of a real symmetric matrix depend continuously on its entries, the uniform continuity of $D^2u$ in time and the uniform positive initial gap imply that $1+\lambda_j\lambda_k>0$ for all $j\ne k$ on a sufficiently short time interval.

By Lemmas~\ref{lem:evolution} and~\ref{lem:cutoff}, $\omega[u](x,t)\ge\inf\omega[u_0]$ and $\Lambda[u](x,t)\ge\inf\Lambda[u_0]>0.$
The first inequality bounds all eigenvalues since
\begin{equation}\label{eq:hessian-from-omega}
 \prod_j(1+\lambda_j^2)\le\omega_0^{-2}
 \,\text{ and }\,
 \sum_j\lambda_j^2\le\omega_0^{-2}-1.
\end{equation}
For the second assertion, we use that each factor in the positive product $\Lambda$ is at most one. 
Thus, each factor is at least $\inf\Lambda[u_0]$, so $1+\lambda_j\lambda_k$ retains a uniform positive lower bound. 

For completeness, rescale about $(x_0,t_0)$ by
\[
 \widehat u(y,s)=t_0^{-1}u(x_0+\sqrt{t_0}\,y,t_0(1+s))
\]
for $s\in[-1,0].$
The Hessian bound is unchanged. 
The interior estimate gives
\[
 |D^3u(x_0,t_0)|\le C(n,M)t_0^{-1/2}.
\]
Applying the differentiated equation and interior Schauder estimates on successively smaller cylinders yields
$|D^\ell u(x_0,t_0)|\le C_\ell(n,M)t_0^{1-\ell/2}$ for $\ell\ge3$.
One may subtract an affine function and a constant on each rescaled cylinder:
the Hessian bound and $|u_t|\le n\pi/2$ then bound its zeroth and first
spatial derivatives there. 
Thus, the Schauder constants are independent of the location of the cylinder and of the initial higher derivatives.

If the maximal time $T$ were finite, these bounds on $[T/2,T)$, together
with the short-time bounds, would bound all derivatives of order at least
two up to $T$. 
The short-time existence theorem would extend the solution past $T$. 
Therefore, $T=\infty$.
\end{proof}

In the general case (Corollary~\ref{cor:preserve-AI}), the uniqueness part is proven based on \cite{CY07}.
Here, we present a short argument that works directly for the graphical case.

\begin{lemma}
	\label{lem:uniqueness}
Let $u$ and $v$ solve \eqref{eq:flow}, be continuous up to $t=0$, satisfy $u(x,0)=v(x, 0)$ for every $x\in \R^n$, and have uniformly bounded Hessians on finite time intervals.
Then $u=v$.
\end{lemma}

\begin{proof}
Since $|u_t|,|v_t|\le n\pi/2$, their difference $w=u-v$ satisfies
$|w|\le n\pi T$ on $\R^n\times[0,T]$. 
The mean value formula applying to $w_t=\Theta(D^2u)-\Theta(D^2v)$ gives $w_t=a^{jk}w_{jk}$ with
\[
 a^{jk}=\int_0^1\big[I+(D^2v+\tau(D^2u-D^2v))^2\big]^{-1}_{jk}\,d\tau.
\]
For a fixed Hessian bound $K$,
$(1+K^2)^{-1}I\le a\le I$. 
Apply the maximum principle to the bounded function $w$, or compare it with
$\eps(1+|x|^2+2nt)$ and then let $\eps\searrow 0$. 
This gives $w\le0$, and interchanging $u$ and $v$ implies the equality. 
The same proof applies to a $C^{1,1}$ common initial value.
\end{proof}

The graphical long-time existence theorem now follows from these smooth cases and an approximation lemma given in Appendix~\ref{sec:approximation}.

\begin{proof}[Proof of Theorem~\ref{thm:existence}]
Apply Lemma~\ref{lem:approx-ai} with $M=\sqrt c$, and denote the smooth approximations by $u_{0,k}$ and their solutions by $u_k$.
Proposition~\ref{prop:smooth} supplies global solutions. 
Their initial Hessians are bounded by $M+1$, so \eqref{eq:hessian-from-omega} and the interior estimates give common bounds depending only on $n$ and $c$.
In particular, $D^2u_k$ is uniformly bounded, and
\[
 |u_k(x,t)-u_{0,k}(x)|\le(n\pi/2)t.
\]
On every compact space-time set, these bounds control $u_k$ and, by spatial
interpolation, $Du_k$. 
The equation gives equicontinuity in time for $u_k$; the Hessian bound gives equicontinuity of the spatial gradients on compact sets as well. 
The positive-time higher derivative bounds and the equation therefore give, after taking a diagonal subsequence,
\[
 u_k\longrightarrow u\quad\hbox{locally uniformly up to }t=0,\text{ and}
 \quad
 u_k\longrightarrow u\quad\hbox{locally smoothly for }t>0.
\]
The first convergence uses $u_{0,k}\to u_0$ locally uniformly and the common time-displacement bound. 
The limit solves \eqref{eq:flow}, satisfies \eqref{eq:time-displacement}, and inherits weak almost isoclinicity and \eqref{eq:derivative-decay}. 
Lemma~\ref{lem:uniqueness} proves uniqueness.

We next recover the exact initial density bound, rather than just the weaker bound provided by the approximations. 
Fix $(x,t)$ and a sufficiently large
$R$. 
Apply Corollary~\ref{cor:strictness} to $\omega[u_k]$ and pass to the limit, using $C^2$ convergence of the initial potentials on $\overline B_R$.
This gives
\[
 \omega[u](x,t)\ge\omega_0
 \left(\frac{R^2-|x|^2-|Du(x,t)|^2-2nt}{R^2}\right)^n.
\]
Letting $R\to\infty$ proves $\omega[u]\ge\omega_0$. 
Notice
$\omega_0\ge(1+c)^{-n/2}$, so all constants depend only on $n$ and~$c$.

Finally, suppose that $u_0$ is strictly almost isoclinic at every point. 
Then $m_R=\min_{\overline B_R}\Lambda[u_0]>0$ for each finite $R$.
Corollary~\ref{cor:strictness} passes through the same approximation to give \eqref{eq:local-strictness} for the limit. 
Choose $R$ so large that the right hand side of \eqref{eq:local-strictness} is positive at the prescribed point. 
Therefore, $\Lambda[u](x,t)>0$.
The quantities $1+\lambda_j\lambda_k$ in the limit are already non-negative, so their strictly positive product makes each of them positive.
\end{proof}

\subsubsection{The bounded-gradient rescaled limit}

Next, we prove the existence of a rescaled limit under a bounded gradient assumption.

\begin{proof}[Proof of Theorem~\ref{thm:rescaled}]
For $p=Du,$ Differentiation of \eqref{eq:flow} yields
\begin{equation}\label{eq:gradient-equation}
 \partial_t p^\alpha=g^{jk}\partial_{jk}p^\alpha
 \,\,\text{ and }\,\,
 (g^{jk})=(I+(D^2u)^2)^{-1}.
\end{equation}
Consequently,
\[
 (\partial_t-g^{jk}\partial_{jk})|p|^2
 =-2g^{jk}\ip{\partial_j p}{\partial_k p}\le0.
\]
The Hessian bound and \eqref{eq:time-displacement} imply that $p$ has at most linear growth on finite time intervals. 
The maximum principle thus applies with a small positive multiple of $e^{Ct}(1+|x|^2)^2$ as a barrier at infinity, giving $|p|\le L$.
The gradients $Du(\cdot,t)$ converge locally uniformly to $Du_0$ as $t\searrow 0$ by the local uniform convergence of the potentials and the uniform Hessian bound.

We next estimate the solution $u$ itself, including its additive constant.
For $\eps>0$, set
\[
 b_\eps(x,t)=L\sqrt{|x|^2+2nt+\eps}.
\]
Its Hessian is non-negative, and
\[
 \partial_t b_\eps=\frac{nL}{\sqrt{|x|^2+2nt+\eps}}
 \ge\tr D^2b_\eps\ge\Theta(D^2b_\eps).
\]
Since $\Theta$ is odd, $-b_\eps$ is a subsolution. 
The Lipschitz initial
data satisfy
$u_0(0)-b_\eps(x,0)\le u_0(x)\le u_0(0)+b_\eps(x,0)$.
Comparison followed by $\eps\searrow 0$ gives
\begin{align*}
 |u(x,t)-u_0(0)|\le L\sqrt{|x|^2+2nt}.
\end{align*}
Note that the comparison is valid on the whole space, since by \eqref{eq:time-displacement}, the excess over either barrier is bounded on each finite time interval, and a small quadratic barrier removes the boundary at infinity.
It follows that
\[
 |u_\mu(x,t)|\le\mu^{-2}|u_0(0)|+
                    L\mu^{-1}\sqrt{|x|^2+2nt}.
\]
Thus, $u_\mu\to0$ locally uniformly. 
All rescaled Hessian and higher derivative bounds are uniform on compact subsets with positive time.
The equation controls mixed time derivatives there. 
Smooth compactness therefore upgrades this convergence to $C^\infty_{\mathrm{loc}}$ for the whole family. 
Formula \eqref{eq:normalized-height} is immediate from $|Du|\le L$. 
Finally, taking $\mu=\sqrt{2t+1}$ and evaluating $u_\mu$ at time
$t/(2t+1)\to1/2$ implies the stated locally smooth convergence of the
normalized graphs.
\end{proof}

\subsubsection{Convergence to the self-expander determined by the blow-down}

Finally, we prove the convergence to an expander under the asymptotic conicality assumption.

\begin{proof}[Proof of Theorem~\ref{thm:expander}]
Let $u$ be the solution supplied by Theorem~\ref{thm:existence}. For
$\mu\ge1$, define
\[
 u_{0,\mu}(x)=\mu^{-2}u_0(\mu x)
 \,\,\text{ and }\,\,
 u_\mu(x,t)=\mu^{-2}u(\mu x,\mu^2t),
\]
a rescaling of $u$.
It is distinct from the approximating solutions $u_k$ used in the proof of Theorem~\ref{thm:existence}.

We first examine the initial potentials. 
Direct differentiation gives
\[
 Du_{0,\mu}(x)=\mu^{-1}Du_0(\mu x)
 \,\,\text{ and }\,\,
 D^2u_{0,\mu}(x)=D^2u_0(\mu x).
\]
In particular, their Hessians have a common bound. 
Moreover,
\[
 u_{0,\mu}(0)=\mu^{-2}u_0(0)
 \,\,\text{ and }\,\,
 Du_{0,\mu}(0)=\mu^{-1}Du_0(0),
\]
so their values and gradients at the origin are bounded uniformly for $\mu\ge1$. 
Taylor's formula therefore bounds both $u_{0,\mu}$ and $Du_{0,\mu}$ on every fixed ball, while the Hessian bound makes the gradients uniformly Lipschitz. 
The Arzel\`a--Ascoli theorem and a diagonal argument show that the family is precompact in $C^1_{\mathrm{loc}}(\R^n)$.
The assumed pointwise convergence to $U_0$ identifies every subsequential limit. 
Consequently,
\[
 u_{0,\mu}\longrightarrow U_0
 \quad\hbox{in }C^1_{\mathrm{loc}}(\R^n)
 \quad\hbox{as }\mu\to\infty.
\]
The common Lipschitz bound for the gradients passes to $DU_0$, giving $U_0\in C^{1,1}(\R^n)$. 
For every fixed $\rho>0$, changing the scaling parameter from $\mu$ to $\mu\rho$ gives
\[
 \begin{split}
 U_0(\rho x)
 &=\lim_{\mu\to\infty}\mu^{-2}u_0(\mu\rho x)
 =\rho^2\lim_{\mu\to\infty}(\mu\rho)^{-2}u_0((\mu\rho)x)
 =\rho^2U_0(x).
 \end{split}
\]
Here we use the existence of the full limit as $\mu\to\infty$.
Thus, $U_0$ is homogeneous of degree two.

Next we establish compactness of the rescaled solutions for positive time.
The chain rule gives
\[
 \partial_tu_\mu(x,t)=u_t(\mu x,\mu^2t)
 \,\,\text{ and }\,\,
 D^2u_\mu(x,t)=D^2u(\mu x,\mu^2t).
\]
Therefore, $u_\mu$ solves \eqref{eq:flow} with initial value $u_{0,\mu}$.
For every integer $\ell\ge2$, \eqref{eq:derivative-decay} gives
\[
 \begin{split}
 |D^\ell u_\mu(x,t)|^2
 &=\mu^{2\ell-4}|D^\ell u(\mu x,\mu^2t)|^2
 \le \mu^{2\ell-4}C_\ell(\mu^2t)^{2-\ell}
 =C_\ell t^{2-\ell}.
 \end{split}
\]
These estimates are independent of $\mu$. 
Rescaling \eqref{eq:time-displacement} also gives
\begin{equation}\label{eq:scaled-time-trace}
 |u_\mu(x,t)-u_{0,\mu}(x)|\le\frac{n\pi}{2}t.
\end{equation}
Since $u_{0,\mu}$ converges locally uniformly, this inequality bounds $u_\mu$ uniformly on every set $B_R\times[0,T]$. 
The common Hessian bound then also bounds $Du_\mu$ there. 
More explicitly, if $\|D^2u_\mu\|_{\op}\le K$, Taylor's formula along any unit vector $e$ gives
\[
 |D_eu_\mu(x,t)|
 \le |u_\mu(x+e,t)-u_\mu(x,t)|+\frac K2,
\]
and the right hand side is uniformly bounded for $(x,t)\in B_R\times[0,T]$ by the bound for $u_\mu$ on $B_{R+1}\times[0,T]$.

On a compact cylinder $B_R\times[\tau,T]$, with $\tau>0$, all spatial derivatives are therefore uniformly bounded. 
Differentiating $\partial_tu_\mu=\Theta(D^2u_\mu)$ supplies uniform bounds for the mixed space and time derivatives as well. 
Smooth compactness and a diagonal argument now show that every sequence $\mu_j\to\infty$ has a subsequence converging in $C^\infty_{\mathrm{loc}}(\R^n\times(0,\infty))$ to a solution $U$ of \eqref{eq:flow}.

Passing to the limit in \eqref{eq:scaled-time-trace}, using
$u_{0,\mu_j}\to U_0$, gives
\[
 |U(x,t)-U_0(x)|\le\frac{n\pi}{2}t
 \qquad\hbox{for every }x\in\R^n,\ t>0.
\]
In particular, $U(\cdot,t)\longrightarrow U_0$ locally uniformly as $t\searrow 0.$
We extend $U$ continuously to $t=0$ by setting $U(x,0)=U_0(x)$.
The Hessian bound and the weak almost isoclinic inequalities pass to the smooth limit for positive time.

We next show that the limit does not depend on the chosen subsequence.
Let $U$ and $\widehat U$ be two possible limits. 
The preceding bound yields
\[
 |U(x,t)-\widehat U(x,t)|\le n\pi t.
\]
Hence, their difference is bounded on $\R^n\times[0,T]$ for each finite
$T$, even though the potentials themselves need not be bounded in space.
Both have uniformly bounded Hessians and converge to $U_0$ locally uniformly as $t\searrow 0$. 
Lemma~\ref{lem:uniqueness} therefore gives $U=\widehat U$.

To spell out this application, put $w=U-\widehat U$. 
For $t>0$, the mean value formula gives $w_t=a^{jk}w_{jk}$ with
\[
 a^{jk}=\int_0^1
 \left[I+\bigl(D^2\widehat U
       +s(D^2U-D^2\widehat U)\bigr)^2\right]^{-1}_{jk}\,ds.
\]
A common operator norm bound $K$ for the Hessians implies
\[
 (1+K^2)^{-1}I\le(a^{jk})\le I.
\]
Fix $0<\tau<T$. 
At time $\tau,$ we have $|w(\cdot,\tau)|\le n\pi\tau$. 
The maximum principle for bounded solutions on $\R^n\times[\tau,T]$ gives
$|w|\le n\pi\tau$ on this interval. 
One can justify the use of the maximum principle on the unbounded domain by comparing $w$ and $-w$ with
\[
 n\pi\tau+\eta\bigl(1+|x|^2+2n(t-\tau)\bigr)
 \,\,\text{ where }\,\,
  \eta>0,
\]
which is a supersolution because
\[
 (\partial_t-a^{jk}\partial_{jk})
 \bigl(1+|x|^2+2n(t-\tau)\bigr)
 =2n-2\operatorname{tr}(a^{jk})\ge0.
\]
Letting $\eta\searrow 0$ and then $\tau\searrow 0$ proves $w=0$.

It follows that the whole family $u_\mu$ converges to $U$, rather than only a subsequence. 
Indeed, if the convergence failed in some $C^m$ norm on a compact subset of $\R^n\times(0,\infty)$, there would be a sequence $\mu_j\to\infty$ staying a fixed positive distance from $U$ in that norm.
The compactness that is already proved would give a subsequence converging smoothly to another limit, which uniqueness identifies with $U$.
This is a contradiction.

It remains to prove the self-similarity. 
For a fixed $\rho>0$, set
\[
 U_\rho(x,t)=\rho^{-2}U(\rho x,\rho^2t).
\]
The same scaling calculation shows that $U_\rho$ solves \eqref{eq:flow} and has a uniformly bounded Hessian. 
The homogeneity of $U_0$ gives
\[
 \begin{split}
 |U_\rho(x,t)-U_0(x)|
 &=\rho^{-2}|U(\rho x,\rho^2t)-U_0(\rho x)|
 \le\frac{n\pi}{2}t.
 \end{split}
\]
Thus, $U_\rho(\cdot,t)\to U_0$ locally uniformly as $t\searrow 0$.
Uniqueness applies once again and yields $U_\rho=U$. 
Equivalently,
\(
 U(\rho x,\rho^2t)=\rho^2U(x,t)
\)
for every $\rho>0.$
Writing $V=U(\cdot,1)$ and taking $\rho=\sqrt t$ gives
\[
 U(x,t)=tV(x/\sqrt t),
\]
which proves \eqref{eq:selfsimilar}.

Finally, with $y=x/\sqrt t$, direct differentiation gives
\[
 U_t(x,t)=V(y)-\frac12y\cdot DV(y)
 \,\,\text{ and }
 D_x^2U(x,t)=D^2V(y).
\]
Substitution into \eqref{eq:flow}, followed by evaluation at $t=1$, proves \eqref{eq:expander-equation}. 
Moreover,
\[
 DU(x,t)=\sqrt t\,DV(x/\sqrt t).
\]
Consequently, after writing $x=\sqrt t\,y$, the gradient graphs satisfy
\[
 \begin{split}
 L_t^U
 &=\{(x,DU(x,t)):x\in\R^n\}
 =\{\sqrt t\,(y,DV(y)):y\in\R^n\}
 =\sqrt t\,L_1^U.
 \end{split}
\]
For the parametrization $F_t(y)=\sqrt t\,(y,DV(y))$, its velocity at $t=1$ is $\tfrac12F_1(y)$. 
Since these graphs evolve by mean curvature up to tangential reparametrization, their normal velocity equals their mean curvature vector. 
Therefore,
\[
\h_{L_1^U}=\frac12F_1^\perp,
\]
which is the self-expander equation.
\end{proof}

\section{\bf Even-dimensional case}
\label{sec:even-G}

\subsection{Graphicality over complex planes}

In this section, we show that any almost isoclinic Lagrangian plane in an even-dimensional complex Euclidean space is graphical over a fixed complex plane. 
More precisely, let \(n=2m\), where \(m\) is a positive integer. 
We show in Theorem~\ref{thm:complex-graphicality} that every almost isoclinic real \(2m\)-dimensional Lagrangian plane \(L\subset\mathbb C^{2m}\) is graphical over the fixed complex \(m\)-plane whose Euclidean volume form is
\begin{equation}\label{eq:Omega}
\Omega=\left(\frac i4\right)^m
\bigwedge_{a=1}^{m}
\pr{(dz_{2a-1}-i\,dz_{2a})
\wedge(d\bar z_{2a-1}+i\,d\bar z_{2a})}.
\end{equation}
Throughout the section, the factors are ordered using the fixed coordinate pairs
\[
(1,2),\ (3,4),\ \ldots,\ (2m-1,2m).
\]
The complex \(m\)-plane in the theorem is
\begin{equation}\label{eq:complex-plane}
P=\{\eta_1=\cdots=\eta_m=0\}
\quad\text{ where }\quad
\eta_a=\frac{z_{2a-1}+i z_{2a}}{\sqrt2}\text{ for } a=1,\ldots,m.
\end{equation}
Indeed, the covectors \((dz_{2a-1}-i\,dz_{2a})/\sqrt2\) give orthonormal complex coordinate covectors on \(P\) and annihilate \(P^\perp\). 
Thus, with its complex orientation, \(\Omega\) is the Euclidean volume form of \(P\), pulled back by orthogonal projection.

We first express $\Omega$ in terms of the standard K\"ahler form and a real 2-form; see \eqref{eq:ambient-volume}.
When \(m=1\), the underlying identity is
\begin{equation*}
\begin{aligned}
2\Omega
&=\frac i2(dz_1-i\,dz_2)\wedge(d\bar z_1+i\,d\bar z_2)
=\frac i2\bigl(dz_1\wedge d\bar z_1+dz_2\wedge d\bar z_2\bigr)
-\frac12\bigl(dz_1\wedge d\bar z_2+d\bar z_1\wedge d z_2\bigr).
\end{aligned}
\end{equation*}
With respect to the standard complex orientation of \(\mathbb C^2\), the first summand on the right is self-dual and the second is anti-self-dual; both have norm \(\sqrt2\).
Consequently, \(\Omega\) is the sum of a self-dual and an anti-self-dual two-form, each of norm \(1/\sqrt2\), and is the Euclidean volume form of a complex one-plane. 
This is the essential linear-algebraic fact used in \cite{W01-MRL}.

To generalize this construction from \(\mathbb C^2\) to \(\mathbb C^{2m}\), we consider the standard K\"ahler form
\begin{equation*}
\omega=\frac i2\sum_{j=1}^{2m}dz_j\wedge d\bar z_j
\end{equation*}
and the real two-form
\begin{equation*}\label{eq:beta}
\begin{aligned}
\beta
&=\frac12\sum_{a=1}^{m}
\bigl(dz_{2a-1}\wedge d\bar z_{2a}
+d\bar z_{2a-1}\wedge d z_{2a}\bigr).
\end{aligned}
\end{equation*}
The form \(\beta\) is the K\"ahler form of a second parallel orthogonal complex structure. 
For each such pair, direct expansion gives
\begin{equation*}
\begin{aligned}
&\frac i2(dz_{2a-1}-i\,dz_{2a})\wedge(d\bar z_{2a-1}+i\,d\bar z_{2a})\\
=&\,\frac i2\bigl(dz_{2a-1}\wedge d\bar z_{2a-1}
+dz_{2a}\wedge d\bar z_{2a}\bigr)
-\frac12\bigl(dz_{2a-1}\wedge d\bar z_{2a}
+d\bar z_{2a-1}\wedge d z_{2a}\bigr).
\end{aligned}
\end{equation*}
Summing over the fixed pairs gives
\begin{equation*}
\omega-\beta=\frac i2\sum_{a=1}^{m}
(dz_{2a-1}-i\,dz_{2a})\wedge(d\bar z_{2a-1}+i\,d\bar z_{2a}).
\end{equation*}
Each summand is decomposable and has square zero, while two-forms commute under wedge product. 
Taking the \(m\)-th exterior power and comparing with \eqref{eq:Omega}, we obtain the identity
\begin{equation}\label{eq:ambient-volume}
\Omega=\frac{(\omega-\beta)^{\wedge m}}{2^m m!}.
\end{equation}

For an oriented real \(2m\)-plane \(L\), we write
\(
*_L\Omega=\Omega(e_1,\ldots,e_{2m}),
\)
where \(e_1,\ldots,e_{2m}\) is any positively oriented orthonormal basis of \(L\). 
Equivalently, \(*_L\) denotes the Hodge star of the pullback to \(L\). 
The number \(*_L\Omega\) is the signed Jacobian of the orthogonal projection \(\pi_P\colon L\to P\). 
If \(L\) is Lagrangian, then \(\omega\) vanishes on pairs of vectors in \(L\). 
Therefore, \eqref{eq:ambient-volume} implies
\begin{equation}\label{eq:star-volume}
*_L\Omega=\frac{(-1)^m}{2^m m!}*_L(\beta^{\wedge m}).
\end{equation}

By Lemma \ref{lem:angle-basis}, we choose an oriented orthonormal basis \eqref{eq:adapted-frame} for $L$ of the form
$e_j=e^{i\theta_j}u_j$ for $j=1,\ldots,n$. 
Write
\[
u_{kj}=dz_k(u_j)\in\mathbb R
\quad\text{ and }\quad
O=(u_{kj})_{k,j=1}^{2m}\in SO(2m).
\]
Then \(dz_k(e_p)=e^{i\theta_p}u_{kp}\) and
\(d\bar z_k(e_p)=e^{-i\theta_p}u_{kp}\), so
\begin{equation}\label{eq:coordinate-evaluation}
\begin{aligned}
\RePart(dz_k\wedge d\bar z_\ell)(e_p,e_q)
&=\RePart\bigl(
e^{i(\theta_p-\theta_q)}u_{kp}u_{\ell q}
-e^{i(\theta_q-\theta_p)}u_{kq}u_{\ell p}\bigr)\\
&=\cos(\theta_p-\theta_q)
\bigl(u_{kp}u_{\ell q}-u_{kq}u_{\ell p}\bigr).
\end{aligned}
\end{equation}
Define the real skew-symmetric matrices
$K=(K_{pq})$ and $B=(B_{pq})$ where
\begin{equation*}
K_{pq}=\beta(u_p,u_q)
\quad\text{ and }\quad
B_{pq}=\beta(e_p,e_q).
\end{equation*}
Summing \eqref{eq:coordinate-evaluation} over the fixed coordinate pairs \((k,\ell)=(1,2),(3,4),\ldots,(2m-1,2m)\) gives
\begin{equation*}
B_{pq}=K_{pq}\cos(\theta_p-\theta_q).
\end{equation*}

For a skew-symmetric \(2m\times2m\) matrix \(A\), recall its Pfaffian is
\begin{equation*}
\Pf A=\frac1{2^m m!}
\sum_{\sigma\in S_{2m}}\operatorname{sgn}(\sigma)
\prod_{a=1}^{m}A_{\sigma(2a-1),\sigma(2a)}.
\end{equation*}
Equivalently, if
\(\alpha_A=\sum_{p<q}A_{pq}\,e^p\wedge e^q\), then
\(
\frac{\alpha_A^{\wedge m}}{m!}
=\Pf A\,e^1\wedge\cdots\wedge e^{2m}.
\)
Applied to \(B\), this gives
\begin{equation}\label{eq:star-Pf}
\frac1{m!}*_L(\beta^{\wedge m})=\Pf B.
\end{equation}
We will also use the identities
\begin{equation}\label{eq:Pf-properties}
(\Pf A)^2=\det A
\quad\text{ and }\quad
\Pf(S^TAS)=\det S\,\Pf A.
\end{equation}

On real vectors, \(dz_j\) and \(d\bar z_j\) have the same values.
Thus, the coefficient matrix of \(\beta\) in the standard real basis is
\begin{equation*}
K_0=\diag\left(
\begin{pmatrix}0&1\\-1&0\end{pmatrix},\ldots,
\begin{pmatrix}0&1\\-1&0\end{pmatrix}\right).
\end{equation*}
The blocks again correspond to the fixed coordinate pairs.
Changing to the positively oriented real basis \(u_j\), we obtain $K=O^TK_0O.$
Consequently,
$K^T=-K,\, K^TK=I,\, K^2=-I,$ and $\Pf K=1.$

Combining \eqref{eq:star-volume} and \eqref{eq:star-Pf}, we obtain
\begin{equation}\label{eq:projection-Pf}
(-1)^m2^m*_L\Omega=\Pf B
=\Pf\bigl[K_{pq}\cos(\theta_p-\theta_q)\bigr]_{p,q=1}^{2m}.
\end{equation}
In the following theorem, we use \eqref{eq:projection-Pf} to show that for an almost isoclinic $L$, $*_L\Omega$ is bounded from below by the isoclinicity function $\Lambda(L)$.

\begin{thm}\label{thm:complex-graphicality}
For any oriented Lagrangian plane \(L\subset\mathbb C^{2m}\) that is almost isoclinic with \(\Lambda(L)>0\), we have
\begin{equation}\label{eq:projection-bound}
(-1)^m*_L\Omega\ge\frac{\Lambda(L)}{2^m}>0.
\end{equation}
In particular, \(L\) is graphical over the fixed complex \(m\)-plane \(P\) in \eqref{eq:complex-plane}, whose Euclidean volume form is \(\Omega\).
\end{thm}

\begin{proof}
The almost isoclinic condition allows us to choose the angles as in Lemma\ref{lem:short-interval}. 
For the following matrix computation, we add a common constant to the angle values so that they lie in \((0,\pi/2)\).
This algebraic substitution preserves their differences, the matrices \(B,K\), and \(\Lambda(L)\). 
Set
\[
C=\diag(\cos\theta_1,\ldots,\cos\theta_{2m})
\quad\text{ and }\quad
T=\diag(\tan\theta_1,\ldots,\tan\theta_{2m}).
\]
Both \(C\) and \(T\) are positive definite. 
We compute
\[
\begin{aligned}
\bigl[C(K+TKT)C\bigr]_{pq}
&=K_{pq}\cos\theta_p\cos\theta_q
\bigl(1+\tan\theta_p\tan\theta_q\bigr)
=K_{pq}\cos(\theta_p-\theta_q)=B_{pq}.
\end{aligned}
\]
Denoting $\widetilde T=KTK^T$ and using $K^2=-I, K^T=-K$, we have
\begin{equation*}
B=C(K+TKT)C=CK(I+\widetilde T T)C.
\end{equation*}
Taking determinants of the both sides and using \eqref{eq:Pf-properties} and
\(\det K=1\), we obtain
\begin{equation}\label{eq:Pf-determinant}
(\Pf B)^2=(\det C)^2\det(I+\widetilde T T).
\end{equation}

Because \(K\) is orthogonal, \(\widetilde T\) is symmetric positive definite with the same eigenvalues as \(T\). 
Write these eigenvalues in increasing order as
\[
0<t_1\le\cdots\le t_{2m}
\]
where $t_j=\tan\theta_{(j)}$ and \(\theta_{(1)}\le\cdots\le\theta_{(2m)}\) is the increasing list of the angle values. 
We keep the matrix indices fixed throughout.
The increasingly ordered eigenvalues of \(T^{-1}\) are \(1/t_{2m},\ldots,1/t_1\). 
Applying
Lemma~\ref{lem:det-rearrangement} below to \(T^{-1}\) and
\(\widetilde T\) gives
\begin{equation}\label{eq:det-estimate}
\begin{aligned}
\det(I+\widetilde T T)
=\det T\,\det(T^{-1}+\widetilde T)
&\ge\left(\prod_{j=1}^{2m}t_{2m+1-j}\right)
\prod_{j=1}^{2m}\left(\frac1{t_{2m+1-j}}+t_j\right)\\
&=\prod_{j=1}^{2m}(1+t_jt_{2m+1-j})\\
&=\pr{\prod_{j=1}^{m}(1+t_jt_{2m+1-j})}^2.
\end{aligned}
\end{equation}

To determine the sign of \(\Pf B\), consider
\[
B_s=C(K+s^2TKT)C
\]
for $0\le s\le1.$
The factorization above gives
\[
\det B_s=(\det C)^2\det(I+s^2\widetilde T T)>0.
\]
In fact, \(\widetilde T T\) is similar to the symmetric positive-definite matrix \(T^{1/2}\widetilde T T^{1/2}\), so all its eigenvalues are positive. 
Thus, \(\Pf B_s\) never vanishes, while
\(\Pf B_0=\det C\,\Pf K=\det C>0\).
Continuity then implies \(\Pf B>0\).

Taking positive square roots in \eqref{eq:Pf-determinant} and
\eqref{eq:det-estimate}, we obtain
\begin{equation*}
\begin{aligned}
\Pf B
\ge\left(\prod_{j=1}^{2m}\cos\theta_{(j)}\right)
\prod_{j=1}^{m}(1+t_jt_{2m+1-j})
&=\prod_{j=1}^{m}
\pr{\cos\theta_{(j)}\cos\theta_{(2m+1-j)}
+\sin\theta_{(j)}\sin\theta_{(2m+1-j)}}\\
&=\prod_{j=1}^{m}
\cos\bigl(\theta_{(j)}-\theta_{(2m+1-j)}\bigr)
\ge\Lambda(L).
\end{aligned}
\end{equation*}
For the last inequality, the displayed \(m\) cosine factors are among those defining \(\Lambda(L)\), and every remaining factor belongs to \((0,1]\). Equation~\eqref{eq:projection-Pf} now proves
\eqref{eq:projection-bound}.

Finally, \(*_L\Omega\ne0\) means that \(\pi_P\colon L\to P\) has nonzero determinant. 
Since both spaces have real dimension \(2m\), this projection is a linear isomorphism. 
Therefore, \(L\) is the graph of the linear map
\(\pi_{P^\perp}\circ(\pi_P|_L)^{-1}\colon P\to P^\perp\).
\end{proof}

We finish this section with the determinant comparison used in the proof.

\begin{lemma}[Fiedler~\cite{Fiedler1971}]\label{lem:det-rearrangement}
Let \(R\) and \(S\) be real symmetric positive-definite \(N\times N\) matrices with increasingly ordered eigenvalues \(r_1\le\cdots\le r_N\) and \(s_1\le\cdots\le s_N\). 
Then
\begin{equation}\label{eq:det-rearrangement}
\det(R+S)\ge\prod_{j=1}^{N}(r_j+s_j).
\end{equation}
\end{lemma}
We include a short proof for completeness. 
\begin{proof}
Minimize \(\log\det(R+USU^T)\) over the compact group \(O(N)\).
At a minimum, write \(V=USU^T\) and \(A=(R+V)^{-1}\).
For every skew-symmetric matrix \(Z\), differentiating along \(e^{tZ}Ve^{-tZ}\) gives
\[
0=\tr\bigl(A(ZV-VZ)\bigr)=\tr\bigl([V,A]Z\bigr).
\]
The commutator \([V,A]=VA-AV\) is skew-symmetric. 
Taking \(Z=[V,A]\), we obtain
\[
0=\tr([V,A]^2)=-\|[V,A]\|^2,
\]
where \(\|\cdot\|\) is the Frobenius norm. 
Hence \([V,A]=0\), so \(V\) also commutes with \(A^{-1}=R+V\), and therefore with \(R\).
Simultaneous orthogonal diagonalization reduces the minimum to one of the products \(\prod_j(r_j+s_{\sigma(j)})\).
For \(r_i\le r_j\) and \(s_a\le s_b\),
\[
(r_i+s_b)(r_j+s_a)-(r_i+s_a)(r_j+s_b)
=(r_j-r_i)(s_b-s_a)\ge0.
\]
Removing inversions shows that the smallest product pairs the two increasing eigenvalue lists in the same order. 
This proves \eqref{eq:det-rearrangement}.
\end{proof}

\subsection{Long-time existence}

Based on Theorem~\ref{thm:complex-graphicality} and Corollary~\ref{cor:tangent-flow-structure}, we use a blow-up argument and the Brakke regularity to prove the long-time existence part in Theorem~\ref{thm:even-flow}.

\begin{proof}
	[Proof of Theorem~\ref{thm:even-flow}]
    For $n\ge 3,$ we apply Corollary~\ref{cor:preserve-AI} to get the flow $M_t$ for a short time.
    For $n=2,$ we apply \cite{SS24}*{Proposition 2.1} to get the flow $M_t$ for a short time, which is valid since Theorem~\ref{thm:complex-graphicality} and the covering map argument as in \cite{LW26}*{Proposition~6.3} imply that there exists a complex plane $P$ given in \eqref{eq:complex-plane} such that $M$ is graphical over $P.$
    In either case, we let $M_t$ be the Lagrangian mean curvature flow starting from $M$ and we let $T\in (0,\infty]$ be the maximal existence time.
	Again, by Theorem~\ref{thm:complex-graphicality} and the covering map argument as in \cite{LW26}*{Proposition~6.3}, $M_t$ is graphical over $P$ for all $t\in [0,T).$
	Note that by \eqref{eq:projection-bound}, each of $M_t$ is the graph of a function with bounded gradient.
	
	Suppose for a contradiction that $T<\infty.$
	Fix a point $x_0\in\bb R^{2n}$ and let $\ovl M_t$ be a tangent flow of $M_t$ at $(x_0,T).$
	Suppose $\ovl M_t$ is non-trivial (i.e., not an empty flow).
	The graphicality and the uniform gradient bound of $M_t$ for $t<T$ implies that each time slice of $\ovl M_t$ is a multiplicity one Lipschitz graph over $P.$
	Corollary~\ref{cor:tangent-flow-structure} then implies that $\ovl M_t$ can only be a multiplicity one planar flow.
	By the Brakke regularity \cite{B78, W02}, $(x_0,T)$ is a regular point of the flow.
	This implies that the flow can be extended smoothly after time $T,$ contradicting the maximailty of $T.$
	Thus, $T=\infty.$
	This proves the long-time existence of the flow.
	
	The convergence of the flow then follows directly from Corollary~\ref{cor:long-time-existence-imply-conv} and Theorem~\ref{thm:rescaled}.
\end{proof}

\section{\bf Odd-dimensional case}
\label{sec:odd}

\subsection{Strengthened bound for $\log\Lambda$}

In this section, we prove a stronger estimate for the evolution of $\log\Lambda$ along an almost isoclinic mean curvature flow when dimension is at least three.
Based on the improvement, a long-time existence result for a Lagrangian mean curvature flow with an almost isoclinic pinching constant will follow.
We remark that the estimates established in this section work in all dimensions $n\ge 3.$

\begin{prop}
\label{prop:evol-improved}
    Let $n\ge 3$ and $M_t$ be a Lagrangian mean curvature flow of almost isoclinic submanifolds in $\bb R^{2n}$ for $t\in I\sbst\bb R.$
    Then the function $\Lambda$ satisfies
    \begin{align}\label{eq:A-43}
    (\bd_t-\D)\log\Lambda
    \ge \frac 43|A|^2
    \end{align}
    for all $t\in I$ and $x\in M_t.$
\end{prop}

\begin{proof}
    We will prove the bound first for $n=3,$ and then show that the three-dimensional bound implies the bound in all dimensions.
  	We will use the bound
	\begin{align}\label{eq:evol-to-be-improved}
		(\bd_t-\D)\log\Lambda
		\ge& \sum_{\ell=1}^n \sum_{j<k}  
		\pr{h_{jj\ell}-h_{kk\ell}}^2
		+ \sum_{\ell=1}^n \sum_{j<k} 
		\pr{4 +\sum_{m\neq j,k} 
			\frac{2(1+\lambda_j\lambda_k)(1+\lambda_m^2)}
			{(1+\lambda_j\lambda_m)(1+\lambda_k\lambda_m)}}
		h_{jk\ell}^2,
	\end{align}
	which follows from \eqref{eq:full-lambda-evolution} based on \eqref{eq:evol-first-term-CS}.

	Assume $n=3$ and we will prove
    \begin{align}\label{eq:n=3-case}
    \sum_{\ell=1}^3 \sum_{j<k}  
		\pr{h_{jj\ell}-h_{kk\ell}}^2
		+ \sum_{\ell=1}^3 \sum_{j<k} 
		\pr{4 +\sum_{m\neq j,k} 
			\frac{2(1+\lambda_j\lambda_k)(1+\lambda_m^2)}
			{(1+\lambda_j\lambda_m)(1+\lambda_k\lambda_m)}}
		h_{jk\ell}^2
    \ge \frac 43|A|^2.
    \end{align}
	At a point $p,$ after rotating it in a good almost isoclinic chart, we may assume $\lambda_1\le \lambda_2=0\le \lambda_3.$
	It follows that $\lambda_1\lambda_3\le 0$ and the almost isoclinic condition and our choice of the char implies $1+\lambda_1\lambda_3>0.$
	Then
	\begin{align*}
		4 + 2\frac{1+\lambda_1^2}{1+\lambda_1\lambda_3}
		&\ge 6 - \frac{2\lambda_1\lambda_3}{1+\lambda_1\lambda_3},\\
		4+2(1+\lambda_1\lambda_3)
		&= 6+2\lambda_1\lambda_3,\text{ and}\\
		4 + 2\frac{1+\lambda_3^2}{1+\lambda_1\lambda_3}
		&\ge 6 - \frac{2\lambda_1\lambda_3}{1+\lambda_1\lambda_3}.
	\end{align*}
	Thus, combining these estimates with the formula 
	\begin{align*}
	3|A|^2-|\h|^2
	= \sum_{\ell=1}^3 \sum_{j<k}  
	\pr{h_{jj\ell}-h_{kk\ell}}^2
	+ \sum_{\ell=1}^3 \sum_{j<k} 6h_{jk\ell}^2
	\end{align*}
	implies that the left hand side of \eqref{eq:n=3-case} is bounded below by
	\begin{align*}
		3|A|^2-|\h|^2 -  \frac{2\lambda_1\lambda_3}{1+\lambda_1\lambda_3} \sum_\ell\pr{h_{12\ell}^2+h_{23\ell}^2}
		+ 2\lambda_1\lambda_3 \sum_\ell h_{13\ell}^2.
	\end{align*}
	Thus, to prove \eqref{eq:n=3-case}, it suffices to prove
	\begin{align}\label{eq:goal-for-4/3}
		|\h|^2
		\le \frac 53|A|^2 -  \frac{2\lambda_1\lambda_3}{1+\lambda_1\lambda_3} \sum_\ell\pr{h_{12\ell}^2+h_{23\ell}^2}
		+ 2\lambda_1\lambda_3 \sum_\ell h_{13\ell}^2.
	\end{align}
	To bound $|\h|^2,$ note that
	\begin{align*}
		\frac 35 
		+ \frac{1+\lambda_1\lambda_3}{5+3\lambda_1\lambda_3}
		+ \frac 1{5+2\lambda_1\lambda_3}
		= \frac{125 + 125\lambda_1\lambda_3 + 28\lambda_1^2\lambda_3^2}{125 + 125\lambda_1\lambda_3 + 30\lambda_1^2\lambda_3^2}
		\le 1,
	\end{align*}
	which implies
	\begin{align*}
		(h_{111}+h_{122}+h_{133})^2
		&\le \frac 53 h_{111}^2 
		+ \frac{5+3\lambda_1\lambda_3}{1+\lambda_1\lambda_3} h_{122}^2
		+ \pr{5+2\lambda_1\lambda_3} h_{133}^2,\\
		(h_{333}+h_{322}+h_{311})^2
		&\le \frac 53 h_{333}^2 
		+ \frac{5+3\lambda_1\lambda_3}{1+\lambda_1\lambda_3} h_{322}^2
		+ \pr{5+2\lambda_1\lambda_3} h_{311}^2,\text{ and}\\
		(h_{222}+h_{211}+h_{233})^2
		&\le \frac 53 h_{222}^2
		+ 5 h_{211}^2
		+ 5 h_{233}^2.
	\end{align*}
	by the Cauchy--Schwarz inequality.
	Taking the sum of these and writing 
	$\frac{5+3\lambda_1\lambda_3}{1+\lambda_1\lambda_3}
	= 5 - \frac{2\lambda_1\lambda_3}{1+\lambda_1\lambda_3}$
	implies
	\begin{align}\label{eq:H-first-bound}
		|\h|^2
		\le \frac 53 |A|^2 
		- \frac{2\lambda_1\lambda_3}{1+\lambda_1\lambda_3} 
		\pr{h_{122}^2+h_{322}^2}
		+ 2\lambda_1\lambda_3\pr{h_{133}^2+h_{311}^2}.
	\end{align}
	The terms to be bounded in \eqref{eq:goal-for-4/3} are
	\begin{align*}
		\frac 53|A|^2 -  \frac{2\lambda_1\lambda_3}{1+\lambda_1\lambda_3} \sum_\ell\pr{h_{12\ell}^2+h_{23\ell}^2}
		+ 2\lambda_1\lambda_3 \sum_\ell h_{13\ell}^2
		=&\,  \frac 53 |A|^2 
		- \frac{2\lambda_1\lambda_3}{1+\lambda_1\lambda_3} 
		\pr{h_{122}^2+h_{322}^2}
		+ 2\lambda_1\lambda_3\pr{h_{133}^2+h_{311}^2}\\
		& - \frac{2\lambda_1\lambda_3}{1+\lambda_1\lambda_3}\pr{h_{121}^2+h_{123}^2+h_{231}^2+h_{233}^2}
		+ 2\lambda_1\lambda_3\cdot h_{132}^2.
	\end{align*}
	Since $\lambda_1\lambda_3\le 0,$ we can throw away the terms with $h_{121}^2$ and $h_{233}^2$ and conclude
	\begin{align*}
		&\frac 53|A|^2 -  \frac{2\lambda_1\lambda_3}{1+\lambda_1\lambda_3} \sum_\ell\pr{h_{12\ell}^2+h_{23\ell}^2}
		+ 2\lambda_1\lambda_3 \sum_\ell h_{13\ell}^2\\
		\ge&\,  \frac 53 |A|^2 
		- \frac{2\lambda_1\lambda_3}{1+\lambda_1\lambda_3} 
		\pr{h_{122}^2+h_{322}^2}
		+ 2\lambda_1\lambda_3\pr{h_{133}^2+h_{311}^2}
		- \frac{2\lambda_1\lambda_3}{1+\lambda_1\lambda_3}\cdot 2h_{123}^2
		+ 2\lambda_1\lambda_3\cdot h_{123}^2\\
		\ge&\, \frac 53 |A|^2 
		- \frac{2\lambda_1\lambda_3}{1+\lambda_1\lambda_3} 
		\pr{h_{122}^2+h_{322}^2}
		+ 2\lambda_1\lambda_3\pr{h_{133}^2+h_{311}^2}
		-2\lambda_1\lambda_3 h_{123}^2\cdot \frac{1-\lambda_1\lambda_3}{1+\lambda_1\lambda_3}\\
		\ge&\, \frac 53 |A|^2 
		- \frac{2\lambda_1\lambda_3}{1+\lambda_1\lambda_3} 
		\pr{h_{122}^2+h_{322}^2}
		+ 2\lambda_1\lambda_3\pr{h_{133}^2+h_{311}^2}.
	\end{align*}
	Hence, \eqref{eq:goal-for-4/3} follows from \eqref{eq:H-first-bound}. 
    This proves \eqref{eq:n=3-case} for every symmetric cubic tensor $h$, and \eqref{eq:evol-to-be-improved} then implies \eqref{eq:A-43} when $n=3$.

    Next, assume $n\ge 3.$
    Let $\cS^n_3:=\set{S\sbst\set{1,2,\cdots,n}:|S|=3}.$
    For any $\set{j,k,m}\in\cS^n_3,$ write
    \begin{align*}
    P_{jkm}
    = P^S_{jk}
    :=\frac{2(1+\lambda_j\lambda_k)(1+\lambda_m^2)}
    	{(1+\lambda_j\lambda_m)(1+\lambda_k\lambda_m)}.
    \end{align*}
    Then \eqref{eq:n=3-case} in the three-dimensional case implies that for any $S\in\cS^n_3,$
    \begin{align}\label{eq:group-by-three}
    \sum_{\ell\in S}
    \sum_{\substack{j,k\in S\\j<k}}
    \pr{
    \pr{h_{jj\ell}-h_{kk\ell}}^2
    + \pr{4+P^S_{jk}} h_{jk\ell}^2
    }
    \ge \frac 43 \sum_{a,b,c\in S}h_{abc}^2.
    \end{align}
    We break the sum \eqref{eq:evol-to-be-improved} into two parts.
    Given $S\in\cS^n_3,$
    by replacing $h_{jk\ell}$ with $0$ for $\set{j,k,\ell}= S$ and applying \eqref{eq:group-by-three}, it follows
    \begin{align*}
    \sum_{\ell\in S}
    \sum_{\substack{j,k\in S\\j<k}}
    \pr{h_{jj\ell}-h_{kk\ell}}^2
    + \sum_{\substack{j,k\in S\\j<k}}
    \sum_{\ell=j,k}
    \pr{4+P^S_{jk}} h_{jk\ell}^2
    \ge \frac 43 \sum_{\substack{a,b,c\in S\\|\set{a,b,c}|\le 2}}h_{abc}^2.
    \end{align*}
    We sum the preceding inequality over $S\in\cS^n_3$ and divide by $n-2$.
    For fixed $j<k$, the difference square $\pr{h_{jj\ell}-h_{kk\ell}}^2$ appears $n-2$ times if $\ell\in\{j,k\}$, since every triple containing $j$ and $k$ contributes this term. 
    If $\ell\notin\{j,k\}$, it appears only once, from the triple $S=\{j,k,\ell\}$. 
    Thus, after dividing them by $n-2$, the coefficient of each difference square is at most one. 
    Moreover, for $\ell\in\{j,k\}$, the averaged coefficient of $h_{jk\ell}^2$ satisfies
    \begin{align*}
    	4+\frac{1}{n-2}\sum_{m\ne j,k}P_{jkm}
    	\le 4+\sum_{m\ne j,k}P_{jkm},
    \end{align*}
    since $P_{jkm}>0$ and $n\ge3$. 
    We then obtain
    \begin{align}\label{eq:c3-jjk}
    \sum_\ell\sum_{j<k} 
    \pr{h_{jj\ell}-h_{kk\ell}}^2
    + \sum_{j<k}
    \sum_{\ell=j,k}
    \pr{4+\sum_{m\neq j,k} P_{jkm}} h_{jk\ell}^2
    \ge \frac{4}{3(n-2)} \sum_{S\in\cS^n_3} \sum_{\substack{a,b,c\in S\\|\set{a,b,c}|\le 2}}h_{abc}^2.
    \end{align}
    For the part with $|\set{j,k,\ell}|=3,$ which contributes
    \begin{align*}
    \sum_{\substack{j,k,\ell\\j<k,\, \ell\neq j,k}}
    \pr{4+\sum_{m\neq j,k} P_{jkm}} h_{jk\ell}^2
    \end{align*}
    in \eqref{eq:evol-to-be-improved},
    we ignore the $P_{jkm}$ part and count directly
    \begin{align}\label{eq:c3-jkl}
    \sum_{\substack{j,k,\ell\\j<k,\, \ell\neq j,k}}
    4 h_{jk\ell}^2
    \ge 12 \sum_{\set{j,k,\ell}\in \cS^n_3}h_{jk\ell}^2
    \end{align}
    as every $(j,k,\ell)\in\cS^n_3$ appears three times in the sum on the left.
    Combining \eqref{eq:c3-jjk} and \eqref{eq:c3-jkl} with~\eqref{eq:evol-to-be-improved} leads to
    \begin{align*}
    (\bd_t-\D)\log\Lambda
    \ge \frac{4}{3(n-2)} \sum_{S\in\cS^n_3} \sum_{\substack{a,b,c\in S\\|\set{a,b,c}|\le 2}}h_{abc}^2
    + 12 \sum_{\set{j,k,\ell}\in \cS^n_3}h_{jk\ell}^2.
    \end{align*}
    For the first part, there are two types $h_{jjj}$ and $h_{jjk},$ and they appear $\binom{n-1}2$ and $3(n-2)$ times, so
    \begin{align*}
    (\bd_t-\D)\log\Lambda
    &\ge \frac{4}{3(n-2)} \sum_{S\in\cS^n_3} \sum_{\substack{a,b,c\in S\\|\set{a,b,c}|\le 2}}h_{abc}^2
    + 12 \sum_{(j,k,\ell)\in \cS^n_3}h_{jk\ell}^2\\
    &= \frac{4}{3(n-2)}
    \pr{
    \sum_j \binom{n-1}2 h_{jjj}^2
    + \sum_{j\neq k} 3(n-2)h_{jjk}^2
    }
    + 12 \sum_{(j,k,\ell)\in \cS^n_3}h_{jk\ell}^2\\
    &= \frac 43 \cdot\frac{n-1}2 \sum_j h_{jjj}^2
    + \frac 43 \cdot \sum_{j\neq k}3h_{jjk}^2 
    + 2 \sum_{(j,k,\ell)\in \cS^n_3}6h_{jk\ell}^2. 
    \end{align*}
    Note that
    \begin{align*}
    |A|^2
    = \sum_j h_{jjj}^2
    + \sum_{j\neq k}3h_{jjk}^2 
    + \sum_{(j,k,\ell)\in \cS^n_3}6h_{jk\ell}^2,
    \end{align*}
    so it follows that
    \begin{align*}
    (\bd_t-\D)\log\Lambda
    \ge \frac 43|A|^2
    \end{align*}
    since $n\ge 3$ and $4/3\le 2.$
\end{proof}

\subsection{Pointwise curvature decay and long-time existence}

In this section, we use the improved bound in Proposition~\ref{prop:evol-improved} to derive a long-time existence result for almost isoclinic Lagrangian mean curvature flow with a pinching constant on $\Lambda.$
In the following proposition, we let 
\begin{align*}
c_0:= e^{-8/9}.
\end{align*}
This constant is approximately $0.41$ and is strictly less than $1/2.$

\begin{prop}
\label{prop:AF-long-time-existence}
    Let $n\ge 3.$
    Suppose $M$ is a complete, connected, almost isoclinic Lagrangian in $\bb R^{2n}$ with a globally defined Lagrangian angle and suppose $M$ is asymptotic to a union of almost isoclinic Lagrangian $n$-planes with the same Lagrangian angles.
    If there exists $\varepsilon>0$ such that $\Lambda\ge c_0+\varepsilon$ on $M,$ then the Lagrangian mean curvature flow starting from $M$ exists for all time.
\end{prop}

\begin{proof}
	From Corollary~\ref{cor:preserve-AI}, there exists a unique Lagrangian mean curvature flow $M_t$ starting from $M$ for $t\in[0,T)$ where $T>0.$
	We will show that the maximal existence time $T$ is infinite by establishing a global curvature decaying estimate.
	
    From \cite{STW}*{Lemma~3.2},
    \begin{align*}
    (\bd_t-\D)\log\pr{2t|A|^2+n}
    \le 3|A|^2 + \frac 12 |\n \log\pr{2t|A|^2+n}|^2.
    \end{align*}
    From the bound $\Lambda\ge e^{-8/9}+\varepsilon,$
    \begin{align}\label{eq:Lambda-lb}
    \log\Lambda> -\frac 89.
    \end{align}
    From Proposition~\ref{prop:evol-improved}, $(\bd_t-\D)\log\Lambda\ge (4/3)|A|^2.$
    Combining these two evolution inequalities, we get
    \begin{align*}
    &(\bd_t-\D) \log\frac{2t|A|^2+n}{\pr{\log\Lambda+8/9}^2}\\
    \le&\, 3|A|^2 + \frac 12 |\n \log\pr{2t|A|^2+n}|^2
    - 2\pr{\frac 1{\log\Lambda+8/9}(\bd_t-\D)\log\Lambda + |\n\log(\log\Lambda+8/9)|^2}\\
    \le&\, 3|A|^2 + \frac 12 |\n \log\pr{2t|A|^2+n}|^2
    - \frac{8/3}{\log\Lambda+8/9}|A|^2
    - \frac 12 |\n\log(\log\Lambda+8/9)^2|^2\\
    =&\, \pr{3-\frac{8/3}{\log\Lambda+8/9}}|A|^2
    + \frac 12\ppair{
    \n\log\frac{2t|A|^2+n}{\pr{\log\Lambda+8/9}^2},
    \n\log\pr{2t|A|^2+n}\pr{\log\Lambda+8/9}^2
    }.
    \end{align*}
    The bound \eqref{eq:Lambda-lb} implies
    \begin{align*}
    3 - \frac{8/3}{\log\Lambda+8/9}
    \le 3 - \frac{8/3}{0+8/9}
    \le 0,
    \end{align*}
    so we derive
    \begin{align*}
    (\bd_t-\D) \log\frac{2t|A|^2+n}{\log\pr{\Lambda+8/9}^2}
    \le \frac 12\ppair{
    	\n\log\frac{2t|A|^2+n}{\pr{\log\Lambda+8/9}^2},
    	\n\log\pr{2t|A|^2+n}\pr{\log\Lambda+8/9}^2
    }.
    \end{align*}
    Based on the maximum principle \cite{EH89}*{Corollary 1.1}, it follows that for any $t\in[0,T),$
    \begin{align*}
    \log\frac{2t|A|^2+n}{\pr{\log\Lambda+8/9}^2}
    \le \sup_{M_0} \log\frac{n}{\pr{\log\Lambda+8/9}^2}
    =: c_n
    \end{align*}
    on $M_t.$
    This implies the existence of a dimensional constant $c_n'<\infty$ such that
    \begin{align*}
    |A|^2\le \frac{c_n'}t
    \end{align*}
    on $M_t$ for all $t\in (0,T),$ and it follows that $T=\infty;$ see, e.g., \cite{S20}*{Proposition~4.6}.    
\end{proof}

We can prove Theorem~\ref{thm:odd-flow} based on the pointwise curvature decay.

\begin{proof}
[Proof of Theorem~\ref{thm:odd-flow}]
    The long-time existence of the flow is proven in Proposition~\ref{prop:AF-long-time-existence}.
    The convergence part then follows directly from Corollary~\ref{cor:long-time-existence-imply-conv} and Theorem~\ref{thm:rescaled}.
\end{proof}

As a corollary, we derive a topological classification for almost isoclinic Lagrangians with $\Lambda>c_0.$

\begin{cor}
\label{cor:AF-top-result}
    Let $n\ge 3.$
    Suppose $M$ is a complete, connected, almost isoclinic Lagrangian in $\bb R^{2n}$ with a globally defined Lagrangian angle and suppose $M$ is asymptotic to a union of almost isoclinic Lagrangian $n$-planes with the same Lagrangian angles.
    If there exists $\varepsilon>0$ such that $\Lambda\ge c_0+\varepsilon$ on $M,$ then $M$ is diffeomorphic to $\bb R^n.$
\end{cor}

The assumptions at infinity and on $\Lambda$ are necessary as there are examples of almost isoclinic submanifolds diffeomorphic to $\bb R\times\bb S^1$ with no positive lower bound on $\Lambda$ and not being asymptotic to a special Lagrangian plane.

\section{\bf Symmetric case}
\label{sec:sym}

In this section, we will prove Theorem~\ref{thm:symmetric-flow} by establishing a distance comparison estimate.
Recall from \cite{LW26}*{Section~4.2} that a submanifold $M\sbst \bb R^{2n}$ is called $SO(n)$-symmetric if there exist an interval $I\sbst\bb R$ and a curve $\gamma=(x,y)\colon I\to\bb R^2$ such that
\begin{align*}
M = \set{\pr{x(r)\,u_1,
	\cdots ,
	x(r)\,u_n,
	y(r)\,u_1,
	\cdots,
	y(r)\,u_n}: r\in I\text{ and }(u_1,\cdots,u_n)\in \bb S^{n-1}},
\end{align*}
where $u=(u_1,\cdots,u_n)\in \mathbb S^{n-1}$ is a parametrization of the standard $(n-1)$-dimensional sphere.
The curve can be chosen to be 
\begin{align}\label{eq:gamma}
\set{(x,y)\in\bb R^2: (x,0,\cdots,0,y,0\cdots,0)\in M}
\end{align}
when $M$ contains the origin.
If $M$ does not contain the origin, the curve \eqref{eq:gamma} may have two components, in which case any of them can serve as $\gamma.$
An $SO(n)$-symmetric submanifold is necessarily Lagrangian~\cite{N07, S20}.
Moreover, it is proven in~\cite{LW26}*{Corollary 4.12} that an almost isoclinic $SO(n)$-symmetric submanifold is always embedded.
Moreover, if the profile curve $\gamma$ is parametrized by $\rho e^{i\phi}\in\bb C$ for some functions $\rho$ and $\phi$ in $r\in\bb I,$ then \cite{LW26}*{Equation~(4.5)} provides
\begin{align}\label{eq:sym-Lambda}
\Lambda
    = \pr{\frac{\rho'}{\sqrt{\pr{\rho'}^2 + \rho^2\pr{\phi'}^2}}}^{n-1}.
\end{align}
In particular, $\rho'\neq 0$ based on the almost isoclinic condition.

First, we show that for a Lagrangian $M$ considered in Theorem~\ref{thm:symmetric-flow}, it must contain the origin.
In particular, $M$ is diffeomorphic to $\bb R^n.$
Recall that we assume all submanifolds are smooth, complete, and connected in $\bb R^{2n}$ unless otherwise stated.

\begin{lemma}
    Let $n\ge 2$ and $M$ be an $n$-dimensional asymptotically conical, $SO(n)$-symmetric, almost isoclinic Lagrangian in $\bb R^{2n}.$
    If there exists $c>0$ such that $\Lambda\ge c>0,$ then $M$ contains the origin, and hence $M$ can be obtained by rotating a curve passing through the origin. 
\end{lemma}

\begin{proof}
    Consider the profile curve
    $\gamma:=\set{(x,y)\in\bb R^2: (x,0,\cdots,0,y,0\cdots,0)\in M}$ parametrized by its arc length $s\in I\sbst\bb R$ and let
    \begin{align*}
    \ell :=\inf_{(x,y)\in\gamma}\sqrt{x^2+y^2}.
    \end{align*}
    We claim that this infimum is achieved.
    Otherwise, there does not exist $s\in I$ such that $\ell = |\gamma(s)|.$
    Then there exists a sequence $s_j\in I$ such that $|\gamma(s_j)|$ is decreasing and
    \begin{align}\label{eq:gamma-s-j-lim}
    \lim_{j\to\infty} |\gamma(s_j)|=\ell
    \end{align}
    The strict monotonicity from \eqref{eq:sym-Lambda} implies that $s_j$ is monotone, so we may assume it is strictly decreasing.
    On the one hand, since $M$ is complete, it follows that $s_j\to-\infty$ as $j\to\infty.$
    On the other hand, using the polar parametrization $\gamma(s) = \rho(s) e^{i\phi(s)},$ formula \eqref{eq:sym-Lambda}, and the bound $\Lambda\ge c>0$ leads to
    \begin{align*}
    \pr{\rho'}^2 \ge c>0,
    \end{align*}
    which contradicts the finiteness in \eqref{eq:gamma-s-j-lim}.

    Thus, there exists $s_0\in I$ such that $\ell = |\gamma(s_0)|.$
    We claim that $\ell=0,$ which will prove the lemma.
    In fact, if $\ell>0,$ formula~\eqref{eq:sym-Lambda} and the condition $\Lambda>0$ would imply that $|\gamma(s)|$ does not have a local minimum at $s_0,$ which contradicts the minimality of $\ell.$
\end{proof}

Thus, an $SO(n)$-symmetric almost isoclinic submanifold is determined by a curve through the origin in $\bb R^2.$
Typical examples of such submanifolds can be found in \cite{LW26}*{Section~4.2};
see also \cite{LW26}*{Figure~3}.

Next, we state the main result of this section.
The following proposition establishes a distance comparison estimate for an $SO(n)$-symmetric flow that passes through the origin.
Along a mean curvature flow of $SO(n)$-symmetric submanifolds $M_t,$ their profile curves $\gamma$ evolve by
\begin{align}\label{eq:modified-CSF}
\pair{\bd_t\gamma,\N}
= k-(n-1)\frac{\pair{\gamma,\N}}{|\gamma|^2}
\end{align}
where $\N$ is a fixed unit normal on $\gamma;$ see, for example, \cite{S20}*{Equation (23)}.
Note that there is a freedom of tangential parametrizations.

\begin{prop}
\label{prop:distance-comp}
    Let $n\ge 2$ and $M$ be an $n$-dimensional asymptotically flat, $SO(n)$-symmetric, almost isoclinic Lagrangian in $\bb R^{2n}.$
    Suppose there exists $c>0$ such that $\Lambda\ge c>0$ and suppose $M_t$ is the unique Lagrangian mean curvature flow starting from $M$ for $t\in[0, T),$ where $T>0.$
    If we parametrize the profile curve $\gamma$ of $M_t$ by their arc lengths, then there exists a tangential parametrization, denoted by $\gamma(s,t)$ for $s\in\bb R$ and $t\in[0,T),$ such that $\gamma(0,t)$ is the origin for all $t\in[0,T)$ and that the following estimate holds.
    There exists $c=c(M)>0$ such that
    \begin{align*}
    |\gamma(s_1,t)-\gamma(s_2,t)
    \ge c_0\cdot |s_1-s_2|
    \end{align*}
    for all $0\le s_1\le s_2.$
\end{prop}

\begin{proof}
	For any $t\in[0,T),$ we assume $\gamma(0,t)$ is the origin and $s$ is the arc length parameter of $\gamma(\cdot,t).$
    We start from choosing better coordinates and tangential parametrizations.

    At a point $\gamma(s,t)$ with $s>0$ and $t\in[0,T),$ write
    \begin{align*}
    \gamma
    = \rho \cdot e^{i\phi}
    = \rho\cdot e_r,
    \end{align*}
    with $e_r=(\cos\phi,\sin\phi)$ and $e_a=(-\sin\phi,\cos\phi).$ 
    Let $\alpha$ be the angle from $e_r$ to $\gamma'$ so that
    \begin{align*}
    \gamma_s=\cos\alpha\cdot e_r + \sin\alpha\cdot e_a.
    \end{align*}
    Then we have
    \begin{align*}
    \cos\alpha\cdot e_r + \sin\alpha\cdot e_a
    = \gamma_s
    = \rho_s e_r + \rho\phi_s e_a,
    \end{align*}
    so $\rho_s=\cos\alpha$ and $\rho\phi_s=\sin\alpha.$
    Thus,
    \begin{align*}
    k = (\alpha+\phi)_s
    = \alpha_s + \frac{\sin\alpha}{\rho},
    \end{align*}
    and hence to get a distance comparison, it suffices to control $\alpha.$

    The normal vector can be chosen to be
    \begin{align*}
    \N = J\gamma_s
    = -\sin\alpha\cdot e_r + \cos\alpha\cdot e_a.
    \end{align*}
    Since $\gamma_t = \rho_te_r + \rho\phi_t e_a,$
    we fix the tangential parametrization such that $\rho_t=0$ (so $\rho$ is fixed in time), which implies
    \begin{align}\label{eq:CSF-1}
    \pair{\gamma_t,\N}=-\rho_t\sin\alpha + \rho\phi_t\cos\alpha
    = \rho\phi_t\cos\alpha.
    \end{align}
    It also follows that $\gamma(0,t)$ being the origin is consistent with this parametrization.
    The almost isoclinic condition allows us to assume $\rho_s>0$ possibly after replacing $s$ with $-s,$ and from now on we assume $s>0.$
    
    The modified curve shortening flow \eqref{eq:modified-CSF} is then controlled by
    \begin{align}\label{eq:CSF-2}
    k-(n-1)\frac{\pair{\gamma,\N}}{|\gamma|^2}
    = \alpha_s + \frac{\sin\alpha}\rho
    - (n-1)\frac{-\rho\sin\alpha}{\rho^2}
    = \alpha_s + n\frac{\sin\alpha}\rho.
    \end{align}
    Thus, combining \eqref{eq:CSF-1} and \eqref{eq:CSF-2} implies
    \begin{align*}
    \rho\phi_t\cos\alpha
    = \alpha_s + n\frac{\sin\alpha}\rho.
    \end{align*}
    Since $\rho_s>0,$ we may use $\rho$ as a new parameter, which satisfies $\alpha_s=\alpha_\rho \rho_s = \alpha_\rho \cos\alpha.$
    It implies 
    \begin{align*}
    \rho\phi_t\cos\alpha
    = \alpha_\rho \cos\alpha + n\frac{\sin\alpha}\rho,
    \end{align*}
    so
    \begin{align*}
    \phi_t
    = \frac{\alpha_\rho}\rho + n\frac{\tan \alpha}{\rho^2}.
    \end{align*}
    Note that $\rho_s>0$ implies $\alpha\in(-\pi/2,\pi/2)$ for all $t\in[0,T).$
    Combining $\rho_t=0$ and $\tan\alpha = (\rho\phi_s)/\rho_s=\rho\phi_\rho$ leads to
    \begin{align*}
    \sec^2\alpha\cdot \alpha_t
    &= \rho\pr{\frac{\alpha_\rho}\rho + n\frac{\tan \alpha}{\rho^2}}_\rho
    = \rho\pr{
        \frac{\alpha_{\rho\rho}}\rho 
        - \frac{\alpha_\rho}{\rho^2}
        + n\frac{\sec^2\alpha\cdot \alpha_\rho}{\rho^2}
        - n\frac{2\tan\alpha}{\rho^3}
    },
    \end{align*}
    which implies
    \begin{align}\label{alpha-evol}
    \alpha_t
    = \cos^2\alpha\cdot \alpha_{\rho\rho}
    + \frac{n-\cos^2\alpha}\rho \alpha_\rho
    - n\frac{\sin 2\alpha}{\rho^2}.
    \end{align}
    
    Now, we apply a maximum principle argument to $\alpha.$
    Recall that the flow is asymptotically flat, which means $\gamma(\cdot,t)$ is asymptotically flat at each time $t\in[0,T).$ 
    Thus, for such $t,$ $\alpha(\cdot,t)$ has a maximum and a minimum.
    At such a maximum point $\rho,$ $\alpha_\rho=0\ge\alpha_{\rho\rho}.$
    Thus, by \eqref{alpha-evol}, at this $\rho,$ 
    \begin{align*}
    \alpha_t\le - n\frac{\sin 2\alpha}{\rho^2}\le 0
    \end{align*}
    since we may assume $\alpha\in (0,\pi/2),$ and the inequality is strict unless the flow is trivial.
    Similarly, at a point $\rho$ when $\alpha$ has a minimum, $\alpha_\rho=0\le\alpha_{\rho\rho}$ so
    \begin{align*}
    \alpha_t\ge - n\frac{\sin 2\alpha}{\rho^2}\ge 0
    \end{align*}
    since we may assume $\alpha\in (-\pi/2,0).$
    These non-zero maximum and minimum points won't happen when $\rho=0$ since $\alpha=0$ when $\rho=0$ by the symmetry.
    This implies $|\alpha|$ has a bound by its supremum at time $t=0$ as long as the flow is smooth.
    
    We use it to get a distance comparison result and conclude the proof.
    For $0\le s_1<s_2,$
    \begin{align*}
    \rho(s_2,t)-\rho(s_1,t)
    = \int_{s_1}^{s_2}\rho_s\,ds
    = \int_{s_1}^{s_2}\cos\alpha\,ds
    \ge c_0(s_2-s_1)
    \end{align*}
    where $c_0:= \inf\limits_{\gamma(\cdot,0)} \cos\alpha$ is a positive constant that only depends on the initial submanifold $M.$
    This particularly implies the desired bound.
\end{proof}

Combining Proposition~\ref{prop:distance-comp} with Corollary~\ref{cor:tangent-flow-structure} implies the long-time existence of the flow.
We then combine it with the graphicality reduction and the convergence result for graphical mean curvature flows to prove Theorem~\ref{thm:symmetric-flow}.

\begin{proof}
[Proof of Theorem~\ref{thm:symmetric-flow}]
    We may assume $n\ge 3$ (or even $n$ is odd) since the two-dimensional case is already covered in Theorem~\ref{thm:even-flow}.
    
    Based on Corollary~\ref{cor:preserve-AI}, let $M_t$ be the Lagrangian mean curvature flow starting from $M.$
    By Proposition~\ref{prop:distance-comp} and Corollary~\ref{cor:tangent-flow-structure}, any non-empty tangent flow of $M_t$ must be a multiplicity one plane, so the flow does not have any finite time singularity by the Brakke regularity \cite{B78, W02}.
    The long-time convergence then follows from Corollary~\ref{cor:long-time-existence-imply-conv} and Theorem~\ref{thm:rescaled}.
\end{proof}

\appendix
\section{Approximation of the initial potential}
\label{sec:approximation}

The following approximation lemma applies even when the initial almost isoclinic inequalities are non-strict.  
Throughout this appendix, eigenvalues are listed in increasing order, and $\|A\|_{\mathrm{op}}$ denotes the operator norm of a symmetric matrix $A$.

\begin{lemma}
\label{lem:approx-ai}
Let $n\geq2$ and $u_0\in C^2(\mathbb R^n)$ satisfy
\[
 \sup_{\mathbb R^n}\|D^2u_0\|_{\mathrm{op}}\leq M
 \,\,\text{ and }\,\,
 1+\lambda_j(D^2u_0)\lambda_\ell(D^2u_0)\geq0
\]
for all $j\neq \ell.$
There is a sequence $u_{0,k}\in C^\infty(\mathbb R^n)$ with the
following properties.
\begin{enumerate}
\item $u_{0,k}\to u_0$ in $C^2_{\mathrm{loc}}(\mathbb R^n)$;
\item $\sup_{\mathbb R^n}\|D^2u_{0,k}\|_{\mathrm{op}}\leq M+1$;
\item for each $k$ and every integer $m\geq2$,
      $\sup_{\mathbb R^n}|D^m u_{0,k}|<\infty$;
\item for each $k,$ there is $d_k>0$ such that on $\mathbb R^n$, for $j\neq\ell,$ $1+\lambda_j(D^2u_{0,k})\lambda_\ell(D^2u_{0,k})\geq d_k.$
\end{enumerate}
If, in addition, there exists $L<\infty$ such that $\sup_{\mathbb R^n}|Du_0|\leq L$, the sequence may
be chosen so that
\[
 \sup_{\mathbb R^n}|Du_{0,k}|\leq (1+o(1))L
\]
as $k\to\infty.$
\end{lemma}

\begin{proof}[Proof of Lemma~\ref{lem:approx-ai}]
Write $b=u_0(0),$ $a=Du_0(0),$ $v(x)=u_0(x)-b-a\cdot x,$ and $r=|x|.$
Taylor's theorem gives
\begin{equation}
 \label{eq:approx-taylor}
 |Dv(x)|\leq Mr,\qquad |v(x)|\leq \frac M2r^2.
\end{equation}
Choose $s_k\in(0,1)$ with $s_k\nearrow 1$ and put
$\gamma_k=1-s_k^2.$
Then choose $\sigma_k\searrow 0$ such that
\begin{equation}
 \label{eq:approx-smallness}
 4\sigma_k\leq1
 \,\,\text{ and }\,\,
 4\sigma_k\leq\frac{\gamma_k}{8(M+1)}.
\end{equation}
Finally, choose $\eta_k>0$ with
$\eta_k\to 0$ and $3(M+1)\eta_k\leq\sigma_k.$

We first construct a cutoff function whose derivatives are small on the logarithmic radial scale.  
Fix a smooth function $\psi\colon\mathbb R\to[0,1]$ that is one on $(-\infty,0]$ and zero on $[1,\infty)$.  
For a sufficiently large $T_k>0,$ set
\[
 \chi_k(x)=
 \begin{cases}
 1,& |x|\leq k,\\
 \displaystyle\psi\!\left(\frac{\log(|x|/k)}{T_k}\right),& |x|>k.
 \end{cases}
\]
This function is smooth, equals one on $B_k$, and vanishes outside $B_{ke^{T_k}}$.
For $r>k$, direct differentiation gives
\[
 r\chi_k'=\frac{\psi'}{T_k}
 \,\,\text{ and }\,\,
 r^2\chi_k''=\frac{\psi''}{T_k^2}-\frac{\psi'}{T_k}.
\]
By increasing $T_k,$ we therefore arrange that
\begin{equation}
 \label{eq:approx-cutoff}
 |r\chi_k'(r)|+|r^2\chi_k''(r)|\leq\eta_k.
\end{equation}
If we let $q_k=s_k\chi_k\in[0,s_k],$ then
$D^2(s_k\chi_kv)=q_kD^2u_0+E_k$
where
\[
 E_k=s_k\bigl(D\chi_k\otimes Dv+Dv\otimes D\chi_k
                    +vD^2\chi_k\bigr).
\]
Since $\|D^2\chi_k\|_{\mathrm{op}}
\leq\max\{|\chi_k''|,|\chi_k'|/r\}$ by the radial symmetry, \eqref{eq:approx-taylor} and \eqref{eq:approx-cutoff}
imply
\begin{align*}
 \|E_k\|_{\mathrm{op}}
 \leq \frac52M\eta_k\leq3M\eta_k\leq\sigma_k.
\end{align*}
The eigenvalues of the main term $q_kD^2u_0$ satisfy, for all $j\neq \ell$,
\begin{equation}
 \label{eq:approx-scaled-products}
 1+q_k^2\lambda_j\lambda_\ell
 \geq 1-s_k^2 
 = \gamma_k.
\end{equation}
This is immediate when $\lambda_j\lambda_\ell\geq0$; otherwise, $\lambda_j\lambda_\ell\geq-1$ and $q_k\leq s_k$ give the assertion.

Consider $g_k=b+a\cdot x+s_k\chi_kv.$
The Hessian $D^2g_k$ vanishes outside a compact set and is continuous everywhere, and hence is uniformly continuous on $\mathbb R^n$.  
Let $\rho_\varepsilon$ be a standard nonnegative smooth compactly supported
mollifier of integral one.  We may choose $0<\varepsilon_k<1/k$ so small
that $u_{0,k}=g_k*\rho_{\varepsilon_k}$ satisfies
\[
 \sup_{\mathbb R^n}\|D^2u_{0,k}-D^2g_k\|_{\mathrm{op}}
 \leq\sigma_k.
\]
The convolutions are well-defined because $g_k$ has at most quadratic growth.  
For $\ell\geq2$, differentiating the mollifier gives
\[
 D^\ell u_{0,k}=D^2g_k*D^{\ell-2}\rho_{\varepsilon_k},
\]
so every derivative of order at least two is globally bounded and hence (3) follows.

The total difference between $D^2u_{0,k}$ and $q_kD^2u_0$ is at most $\varrho_k=4\sigma_k$ in operator norm. 
Their ordered eigenvalues therefore differ by at most $\varrho_k$.

Since the eigenvalues of $q_kD^2u_0$ have absolute value at most $M$, \eqref{eq:approx-smallness} and \eqref{eq:approx-scaled-products} imply, for all $j\neq \ell$,
\[
 1+\lambda_j(D^2u_{0,k})\lambda_\ell(D^2u_{0,k})
 \geq\gamma_k-2M\varrho_k-\varrho_k^2
 \geq\frac{\gamma_k}{2},
\]
where the last inequality also uses $0<\gamma_k<1$.  
Thus, one may take $d_k=\gamma_k/2$ and obtain (4).

The same estimates as above give $\|D^2u_{0,k}\|_{\mathrm{op}}\leq M+1,$ which proves (2).
On a fixed compact set, $\chi_k=1$ for all sufficiently large $k$.
The facts that $s_k\to1$, $\sigma_k\to0$, and $\varepsilon_k\to0$ then imply $u_{0,k}\to u_0$ in $C^2$ on that compact set.
This proves (1).

Finally, suppose $|Du_0|\leq L$.
Then $|a|\leq L$ and $|v(x)|\leq2Lr$, while
\[
 Dg_k=(1-s_k\chi_k)a+s_k\chi_kDu_0+s_kvD\chi_k.
\]
The first two terms form a convex combination of vectors of norm at most $L$.  
The remaining term has norm at most $2L\eta_k$, so
\[
 \sup_{\mathbb R^n}|Dg_k|\leq L+2L\eta_k.
\]
Non-negative convolution preserves this bound, so this proves the bounded-gradient assertion and completes the proof.
\end{proof}

\end{document}